\documentclass[preprint,3p]{elsarticle}
\usepackage{amsmath,amsthm}
\usepackage{amssymb}
\usepackage{array,ulem}

\usepackage[utf8]{inputenc}
\usepackage[T1]{fontenc}
\usepackage{url}
\usepackage[Algorithm]{algorithm} 
\usepackage[noend]{algorithmic}
\usepackage[dvipsnames]{xcolor}
\usepackage{graphicx}
\usepackage{grffile}
\usepackage{amsfonts}
\usepackage{amssymb}
\usepackage{amsmath}
\usepackage{cancel}
\usepackage{bbm}
\usepackage{color}
\usepackage{animate,tabularx,rotating}
\usepackage{IEEEtrantools}
\usepackage{pgfplots}
\usepackage{pgfplotstable}
\usepgfplotslibrary{groupplots}
\usepgfplotslibrary{external}
\usetikzlibrary{plotmarks}
\pgfplotsset{compat=1.16}
\tikzsetfigurename{diffredist-figure}
\usepackage{subcaption}

\usepackage[pagebackref,pdfencoding=auto, psdextra]{hyperref}

\usepackage{todonotes}
\renewcommand{\todo}[2][]{\tikzexternaldisable\@todo[#1]{#2}\tikzexternalenable}

\hypersetup{
    colorlinks=true,       
    linkcolor=black,          
    citecolor=blue,        
    urlcolor=blue,           
}
\usepackage{cleveref}
\crefname{algorithm}{alg.}{alg.}

\newcommand{\eps}{\varepsilon}

\newcommand{\ze}{\frac1\eps\zeta(\frac\phi\eps)}

\newcommand{\x}{\mathbf{x}}   

\newcommand{\Om}{\Omega}
\newcommand{\Hv}{\mathcal{H}}
\renewcommand{\div}{\operatorname{div}}

\newcommand{\off}[1]{}

\newtheorem{lemma}{Lemma}
\newtheorem{proposition}{Proposition}

\newtheorem{remark}{Remark}

\def\R{\mathbb{R}}

\def\Xint#1{\mathchoice
{\XXint\displaystyle\textstyle{#1}}%
{\XXint\textstyle\scriptstyle{#1}}%
{\XXint\scriptstyle\scriptscriptstyle{#1}}%
{\XXint\scriptscriptstyle\scriptscriptstyle{#1}}%
\!\int}
\def\XXint#1#2#3{{\setbox0=\hbox{$#1{#2#3}{\int}$ }
\vcenter{\hbox{$#2#3$ }}\kern-.6\wd0}}

\def\dashint{\Xint-}

\newcommand{\review}[1]{{\color{red}{#1}}}
\newcommand{\reviewII}[1]{{\color{violet}{#1}}}

\crefformat{footnote}{#2\footnotemark[#1]#3}

\renewcommand{\review}[1]{{#1}}
\renewcommand{\reviewII}[1]{{#1}}

\begin{document}

\title{Diffusion-redistanciation schemes for  2D and 3D constrained Willmore
flow: application to the equilibrium shapes of vesicles}

\author[1]{Thibaut Metivet\corref{cor1}}
\ead{thibaut.metivet@inria.fr}
\author[1]{Arnaud Sengers}
\author[2]{Mourad Isma\"\i l}
\author[1]{Emmanuel Maitre}
\cortext[cor1]{Corresponding author}
\address[1]{Univ. Grenoble Alpes, Inria, CNRS, Grenoble INP, LJK, 38000 Grenoble, France}
\address[2]{Laboratoire Interdisciplinaire de Physique, Univ. Grenoble Alpes and CNRS.}

\date{Received: date / Accepted: date}

\begin{abstract}

In this paper we present a novel algorithm for simulating geometrical flows, and
in particular the Willmore flow, with conservation of volume and area.  The idea
is to adapt the class of diffusion-redistanciation algorithms to the Willmore
flow in both two and three dimensions.  These algorithms\footnote{\label{note:algTypo}A typo in the algorithms has been corrected in the current version of the article.} rely on alternating
diffusions of the signed distance function to the interface and a
redistanciation step, and with careful choice of the applied diffusions, end up
moving the zero level-set of the distance function by some geometrical quantity
\review{without resorting to any} explicit transport equation.  The constraints are
enforced between the diffusion and redistanciation steps via a simple rescaling
method. The energy globally decreases at the end of each global step. The
algorithms feature the computational efficiency of thresholding methods without
requiring any adaptive remeshing thanks to the use of a signed distance function
to describe the interface. This opens their application to dynamic
fluid-structure simulations for large and realistic cases.  The methodology is
validated by computing the equilibrium shapes of two- and three-dimensional
vesicles, as well as the Clifford torus.

\end{abstract}

\maketitle

\section{Introduction}
Many different modelling situations can be cast as mathematical problems where an
interface motion is driven by the \review{minimisation} of a geometric energy under
geometric constraints. This is the case of multiphase flows, image segmentation,
or elastic interface modeling, to cite a few. The geometric quantities involved
in energy and constraints are for instance mean or Gaussian curvature of
interface, surface area or enclosed volume. The algorithms and methods developed
in this article give an efficient framework to address such situations.

\review{Our work was motivated by one of the} modeling situations that we now describe.
Vesicles are systems of two fluids separated by a bi-layer membrane of
phospholipid molecules. These objects can be considered as a simple model for
Red Blood Cells (RBC). Since the number of such molecules is constant, this kind
of interface has constant area. Therefore, its shape is determined by high order
energy, i.e. the mean curvature is minimised. Moreover, there is no exchange
across this interface, so that the enclosed volume is constant.  Mathematically,
the problem of \review{finding a surface minimising its mean curvature} is
the well known
problem of Willmore \cite{simon1993existence}. In this work however, we are
interested in this minimisation under the conservation of area and enclosed
volume.

The numerical simulation of vesicles involve the resolution of two-fluid flows
(for the inner and outer fluids) and fluid-structure interactions (for the
membrane-fluids interaction), which is quite challenging: as a sharp object, a
singularity occurs across the membrane bringing stress jump, which should either
be dealt with explicit jump conditions, or appropriate numerical spreading.
Since the membrane energy is of high geometrical order, its gradient involves
high order derivatives of the unknowns. In addition, the inextensibility of the
\review{immersed interface is usually accounted} for using elastic tension energies with
high \review{modulus which makes} the resulting numerical problem very stiff. 

A lot of such numerical simulations have been carried out by different
teams using many numerical methods. To mention some representative works, we can
cite the dynamic molecular method \cite{markvoort2006vesicle}, the boundary
integral method \cite{beaucourt2004steady}, the phase field method (see for
example
\cite{biben2005phase,wang2008modelling,du2008adaptive,lowengrub2009phase,
bretin2015phase}) or level set method \cite{biben2003tumbling,
lowengrub2007surface, maitre2009applications, salac2011level,
doyeux2013simulation}. \cite{kaoui2011two} also proposes a model for
vesicles implemented using a Lattice Boltzmann Method.

In the context of finite-element methods coupled with level set
technique we can cite
\cite{laadhari2014computing,doyeux2013simulation}.  We can also mention
the work \cite{ismail2014necklace} based on a finite-element method
where the membrane is modelled as a necklace of small rigid particles.

In the work above, the nonlinear coupling between the fluid flow and the
geometric \review{description of the interface} is usually made explicit, which leads to
severe restriction on the time step during simulation. Or it could be solved
implicitly by a Newton-type method, which increases dramatically the cost per
iteration. In \cite{cottet2016semi}, a semi-implicit scheme was proposed where
an ad hoc diffusion equation was used as a predictor step for the future
position of a \review{drop of liquid} or a simple elastic interface. 

Our approach in \review{this work has as} primary aim to propose a systematic way to
built a predictor of the position of an interface with constant area and
enclosed volume, when it moves to minimise high \review{geometrical} order energies, such
as its mean curvature. In order to build such a predictor, we extended the
diffusion-thresholding/redistancing approach introduced in \cite{merriman1992diffusion}
and extended in \cite{grzhibovskis2008convolution, esedoglu2008threshold,
esedog2010diffusion, kublik2011algorithms} to the case of area and enclosed
volume conservation.

\paragraph{Outline} The paper is organised as follows: in section \ref{nshgmi}
\review{we review} diffusion-thresholding and diffusion-redistancing schemes principles,
and introduce a new methodology to write higher order motion schemes, such as
the Willmore flow, \review{that} is more systematic than in
\cite{esedog2010diffusion}.  Then we present a new efficient method to take into
account the area and enclosed volume constraints. In section \ref{areavolcons}
we address the problem of volume and area conservation, introducing an explicit
analytic method to project on the constraint set. Then a section is devoted to
numerical aspects of the implementation within finite-element approaches. We
investigate in detail how to choose an optimal time step and plot numerical
convergence curves for a basic scheme and an enhanced one. To evaluate the
performance of the diffusion-redistanciation scheme without rescaling, we study
the convergence of a torus under unconstrained Willmore flow to the optimal
Clifford torus. \review{To conclude}, we present some numerical illustrations of the
computation of 2D and 3D equilibrium shapes for vesicles that match well with
those obtained with classical numerical schemes.

\section{Numerical schemes for higher geometrical motion of interfaces}
\label{nshgmi}
The class of \textit{diffusion-generated motion} was introduced by the work
\cite{merriman1992diffusion} of Merriman, Bence and Osher. They proposed an
efficient algorithm for computing the mean curvature flow of a surface (a curve
in two dimensions) without any
direct computation of the mean curvature. The algorithm consists in repeating
two steps, namely a diffusion step (also named convolution step) and a
thresholding step. During the diffusion step, the characteristic function
representing the surface is diffused for a certain time step, and the
$\frac{1}{2}$ iso-level moves \review{proportionally to the local curvature}. A characteristic
function is then recovered by thresholding the resulting function at $\frac{1}{2}$, 
allowing to iterate the process.

\begin{algorithm}
    \caption{Original Convolution-Thresholding scheme}
    \begin{algorithmic}
        \WHILE{$t < t_f$}
        \STATE{Solve $\partial_t \phi - \Delta \phi = 0$ with initial condition
            $\phi_i = \mathbbm{1}^{(n)}$ for time $\delta t$}
            \STATE{Construct the new characteristic function $\mathbbm{1}^{(n+1)} = \mathbbm{1}_{\lbrace \phi \geq \frac{1}{2} \rbrace}$.}
        \ENDWHILE
    \end{algorithmic}
\end{algorithm}

The main advantage of this method is its simplicity and unconditional stability:
the convolution step consists in the \review{numerical resolution of a heat equation} which can be
achieved efficiently, while the standard phase field approach would involve
solving a more complicated non-linear equation or computing the curvature and
solving a transport equation.

Due to its conceptual simplicity and numerical efficiency, various extensions of
the method have been proposed since its introduction. 
\review{Extension to the multiphase case is trivial for symmetric junctions as one only has to diffuse each interface and threshold to the biggest one in each region. The case of nonsymmetric junctions and different surface tensions is treated in  \cite{ruuth1998diffusion}. However the resulting algorithm relies on a spatially dependent thresholding that require essentially to locate triple junctions. Also the case where four or more phases meet is only treated empirically. In \cite{esed2015threshold} a reformulation of the problem is done and the associated algorithm changes the diffusion step to take into account the difference between the surface tensions.} An area preserving
motion by mean curvature can be found in \cite{ruuth2003simple} by
changing the threshold value from $\frac{1}{2}$ to a real number $\lambda$ found
using a Newton method. We refer to the next section for more details about area
conservation.

Extensions of the convolution-thresholding method to the Willmore flow and other
higher order geometric motion, such as surface diffusion, have been proposed
independently in \cite{esedoglu2008threshold,grzhibovskis2008convolution} in
both two and three dimensions. We recall that the Willmore flow is given by the
following normal velocity : 
\begin{equation} 
    \label{A2:Eqn:VelocityWillmore} 
    W =
    \left\{ 
        \begin{aligned} 
            & \Delta_\Gamma H + \frac{H^3}{2} &&\textrm{ in 2D}
            \\ 
            & \Delta_\Gamma H + 2 H \left( H^2 - K\right) &&\textrm{ in 3D},
        \end{aligned} 
\right.  
\end{equation}
\review{where 
    \begin{equation*}
        H =
        \left\{ 
            \begin{aligned} 
                & \kappa_1 &&\textrm{ in 2D}
                \\ 
                & \frac{\kappa_1+\kappa_2}{2} &&\textrm{ in 3D},
            \end{aligned} 
        \right.  
        \quad\textrm{and}\quad
        K =
        \left\{ 
            \begin{aligned} 
                & \kappa_1 &&\textrm{ in 2D}
                \\ 
                & \kappa_1\kappa_2 &&\textrm{ in 3D},
            \end{aligned} 
        \right.  
    \end{equation*}
    are respectively the mean and Gaussian curvatures of the surface, $\kappa_i$
    denoting the usual principal curvatures.
}

\review{In 2D \cite{esedoglu2008threshold} and in 3D \cite{grzhibovskis2008convolution}}, a local
expansion of the convoluted characteristic function is performed and shows that
the velocity $W$ appears at the second order of the diffused interface. More
precisely,\review{ in the 2D case,} any point of the interface can be relocated at the origin in a way
that the normal to the curve is aligned with the y-axis. The behavior of the
point of interest during the convolution step is then described by the following
expansion:
\reviewII{
\begin{equation}
    \label{A2:eqn:expansionCharacteristicFunction}
   F(0,y,\delta t) = \frac{1}{2} - \frac{1}{2\sqrt{\pi}} y \delta t^{-\frac{1}{2}} 
    + \frac{1}{2\sqrt{\pi}} \kappa \delta t^{\frac{1}{2}} 
    - \frac{1}{4\sqrt{\pi}}W\delta t^\frac{3}{2} 
    + \mathcal{O}\left( \delta t^\frac{5}{2} \right)
\end{equation}
where $F(0,y,\delta t)$ denotes the function obtained by diffusing the
characteristic function for a time $\delta t$.
Thresholding the result at $F(0,y,\delta t) = \frac{1}{2}$ then results in a
normal displacement of the point by $\kappa \delta t$, which gives the classical 
mean curvature scheme. As the
velocity term of interest for the Willmore flow can be found at the second order
in expansion \cref{A2:eqn:expansionCharacteristicFunction}, one can then compute
this flow by extracting it with a linear combination of two different
solutions of the convolution step taken at different time steps $\sqrt{2\Delta t}/\theta$ and $\theta\sqrt{2\Delta t}$ for $\theta\neq 1$ that eliminates
the curvature term. Upon thresholding at an appropriate level \cite{esedoglu2008threshold} found that the curve moved following the Willmore flow since the solution is   \begin{equation*}
    y = \delta t W  + \mathcal{O}\left( \delta t^\frac{3}{2} \right).
\end{equation*}
}
At the expense of an additional diffusion at each step and an error of order
$\delta t^\frac{3}{2}$, we can write an algorithm for the Willmore flow similar to the
mean curvature one. \review{A similar expansion in the 3D case can be found in \cite{grzhibovskis2008convolution} and yields the same expression to compute the Willmore flow.} 

\begin{algorithm}
    \caption{Convolution-Thresholding scheme for Willmore flow}
    \begin{algorithmic}
        \WHILE{$t < t_f$}
            \STATE{Perform one step of $\partial_t \psi - \Delta \psi = 0$ with initial condition
            $\psi_0 = \mathbbm{1}^{(n)}$ for  time steps $\delta t_1=\sqrt{2\delta t}/\theta$ and $\delta t_2=\theta\sqrt{2\delta t}$}
            \STATE{Compute $\mathcal{D} = (2\delta t)^\frac14[\sqrt\theta\psi\left(\delta t_1\right)-\frac1{\sqrt\theta}\psi\left(\delta t_2\right)] $}

            \STATE{Construct the new characteristic function $\mathbbm{1}^{(n+1)} = \mathbbm{1}_{\lbrace \mathcal{D} \geq \frac{1}{2}(2\delta t)^\frac14\left(\sqrt\theta-\frac1{\sqrt\theta}\right) \rbrace}$}
            \STATE{$t \leftarrow t+\delta t$}
        \ENDWHILE
    \end{algorithmic}
\end{algorithm}

This method seems highly valuable because it allows to compute a fourth-order flow
\review{without any differentiation}, especially \review{in the context of
finite-element methods with low-order polynomial discretisations of the
phase-field function $\phi$}.

However the convolution-thresholding method has some identified drawbacks. Its
major issue is its inaccuracy and its propensity to get ``stuck'' when the mesh is
not refined enough at the interface. More precisely, as pointed out in
\cite{merriman1994motion}, given a fixed uniform grid, if the motion during an
iteration is smaller than the grid resolution, the thresholding step can reset
the interface to its initial configuration and loop indefinitely. A solution to
this problem can be found in \cite{ruuth1998efficient}, where the authors use
adaptive grids to refine the mesh near the interface where the resolution is
more important. Such adaptive strategy is also essential for accuracy
considerations, especially in the case of higher order flows, but can become
\review{algorithmically complicated due to its use of unequally spaced FFTs} in the dynamic case or in dimension 3.

An alternative approach to the problem is to allow the function to contain
subgrid informations. In \cite{esedog2010diffusion}, the authors replace to that
end the characteristic function by a signed distance function and explore the
motions that can be obtained through the diffusion and the redistanciation of a
signed distance function. A similar expansion as
\cref{A2:eqn:expansionCharacteristicFunction} of the solution of the heat
equation  starting from a distance function shows that the first order term is
still a mean curvature term. At the expense of a redistanciation step, though more costly
than the thresholding step, one can obtain significantly more accurate computations
than the classical convolution thresholding algorithm. However, second order
term in this case is not the desired Willmore term \review{any more} and an additional correction
is therefore necessary to compute the Willmore flow.

From a theoretical point of view, convergence of the classical algorithm of
motion by mean curvature is well understood and various proofs of the
consistency of the scheme exist. Some are based on a comparison principle
\cite{evans1993convergence,barles1995simple,chambolle2006convergence}, other are
based on a gradient flow interpretation
\cite{almgren1993curvature,luckhaus1995implicit}. Introduced in the two papers
\cite{esed2015threshold,laux2016convergence}, a new interpretation of the
convolution-thresholding algorithm as a gradient minimising flow has initiated
more recent works on the subject \cite{laux2019thresholding,laux2020brakke},
while \cite{swartz2017convergence} relies on a formal matched asymptotic method.
The extension to the multiphase flow can be found in \cite{laux2016convergence}
and the case of volume-preserving motion is studied in
\cite{lauxswartz2016convergence,mugnai2016global}. The recent results on the
gradient flow approach are reviewed in \cite{laux3gradient}. However, the
convergence of higher order motions, especially the Willmore flow, remains an
open question.

In this section, we propose a new local expansion of the signed distance
function near the interface. It yields similar results to the one in
\cite{esedog2010diffusion} for the 2D case but our approach seems more intrinsic
and its generalisation to higher dimensions and orders is straightforward. We propose
diffusion-redistanciation algorithms to compute the Willmore flow in 2D and
3D. As the term in the second order of the expansion is not quite the Willmore
term, one has to add a correction term. The 3D version in particular shows a
completely new approach and
features a correction term using a diffusion of $d^2$.

\subsection{Expansion of a distance function convoluted by the heat kernel}
\label{subsec:convolutionexpansion}

As suggested by the previous work of Esedoglu, Ruuth and Tsai
\cite{esedog2010diffusion}, the high-order geometrical motions of a surface can be
computed using the convolution of a distance function to this interface by the
heat equation kernel. In this section, we thus derive an expansion of the
solution of the heat equation with a distance function as its initial condition.
In contrast with \cite{esedog2010diffusion}, we compute this expansion for both the
two- and three-dimensional cases using a fully implicit approach.

We consider $d$ a signed distance function to some interface $\Gamma$,
and introduce the heat equation:
\begin{equation}
    \partial_t \phi - \Delta \phi = 0.
    \label{eq:heat_equation}
\end{equation}
supplemented with the initial condition: 
\begin{equation}
    \phi(t=0) = d.
    \label{eq:heat_initial_condition}
\end{equation}

The Taylor expansion of the heat equation solution with respect to time then reads:
\begin{equation}
    \phi( \delta t ) = \phi(0) 
    + \partial_t \phi(0) \: \delta t 
    + \partial_{tt} \phi(0) \: \frac{\delta t^2}{2} 
    + \mathcal{O}\left( \delta t^3 \right)
    \label{eq:Gt_expansion_1}
\end{equation}

Recalling \cref{eq:heat_equation} and the initial condition
\cref{eq:heat_initial_condition}, we then have $\partial_t
\phi(0) = \Delta d$ and $\partial_{tt} \phi(0) = \Delta^2 d$.
We denote $w = \frac{\Delta^2 d}{2}$ and can then rewrite \cref{eq:Gt_expansion_1} as:
\begin{equation}
    \begin{aligned}
        \mathcal{G}_{\delta t} \equiv \phi( \delta t )
        & = d + \Delta d \: \delta t + \frac{\Delta^2 d}{2} \: \delta t^2 + \mathcal{O}\left(
        \delta t^3\right) \\
        & = d + \Delta d \: \delta t + w \: \delta t^2 + \mathcal{O}\left( \delta t^3
        \right).
    \end{aligned}
    \label{eq:Gt_expansion_2}
\end{equation}

\begin{remark}
    \label{rk:rem_time_discret}
    When time discretisation of (\ref{eq:heat_equation}) will come into play, we
    will have to consider a second order scheme to get the right expansion. Using a
    first order Euler scheme, which could be desirable for \review{its} stability properties,
    leads to double the second order term magnitude. We refer to
    \cref{solving_the_diff_eq} for more details.  
\end{remark}

By definition, $\Delta d$ is the total curvature of the local $d$ iso-surface:
\begin{equation*}
    \Delta d = \kappa = (n-1) H
\end{equation*}
with $n$ the space dimension and $H$ the mean curvature.
We can also express $\Delta^2 d$ close to the interface $\Gamma$ as a function
of the interface geometrical quantities, and in particular the
\review{Laplace-Beltrami operator} of the curvature $\Delta_{\Gamma} \kappa = 
(n-1)\Delta_{\Gamma} H$. 
More specifically, \review{if $p$ denotes the projection of $x$ on $\Gamma$, } we have $\Delta_{\Gamma}f = \Delta\left(f(x - d\nabla
d)\right)|_{\Gamma} \equiv \Delta\left(f \circ p\right)|_{\Gamma}$ for any 
$f \in \mathcal{C}^2(\Gamma)$ \cite{delfour2011shapes}, so that:
\begin{equation*}
    \begin{aligned}
        \Delta_{\Gamma}(\Delta d) 
        & = \Delta(\Delta d \circ p)|_{\Gamma} \\
        & = \nabla_i \left( \nabla_j p_i \, \nabla_j \Delta d \right)|_{\Gamma} \\
        & = \nabla_i \left( (\delta_{ij} - \nabla_i d \, \nabla_j d - d\nabla_i
        \nabla_j d) \, \nabla_j \Delta d \right)|_{\Gamma} \\
        & = \left[ -\Delta d \, \nabla_{{N}} \Delta d - \nabla_{{N}}^2
            \Delta d + \Delta^2 d
        - d \nabla_i \left( \nabla_i\nabla_j d \, \nabla_j \Delta d \right)
        \right]|_{\Gamma}
    \end{aligned}
\end{equation*}
where $N \equiv \nabla d$ is the interface normal and we have introduced
the normal gradient $\nabla_{{N}} \equiv {N}\cdot\nabla$. Applying the
restriction to $\Gamma \equiv \left\{d = 0\right\}$, we finally obtain:
\begin{equation}
    \Delta^2 d = \Delta_{\Gamma} \kappa + \nabla_{{N}}^2 \kappa + \kappa \,
    \nabla_{{N}} \kappa
    \label{eq:delta2_d_eq1}
\end{equation}
We can also compute:
\begin{equation*}
    \begin{aligned}
        \nabla_{{N}} \kappa
        & = \nabla_i d \: \nabla_i \nabla_k \nabla_k d \\
        & = \nabla_i d \: \nabla_k \nabla_k \nabla_i d \\
        & = \nabla_k \left( \cancel{\nabla_i d \: \nabla_k \nabla_i d} \right) 
        - \nabla_k \nabla_i d \: \nabla_i\nabla_k d \\
        & = - \mathrm{Tr}\left[ (\nabla^2 d)^2 \right] \\
    \end{aligned}
\end{equation*}
since the distance property $|\nabla d| = 1$ gives $\nabla_i d \: \nabla_k \nabla_i d = \nabla_k |\nabla d|^2 / 2 = 0$.
Recalling that $\nabla^2 d$ is the Weingarten map of $\Gamma$, we can then write:
\begin{equation}
    \nabla_{{N}} \kappa = -\sum_i \kappa_i^2
    = \left\{
    \begin{aligned}
        & -\kappa^2
        \\
        & 2K - \kappa^2
    \end{aligned}
    \right.
    = \left\{
    \begin{aligned}
        & -H^2 &\textrm{ in 2D}
        \\
        & -2 \left(2 H^2 - K\right) &\textrm{ in 3D}
    \end{aligned}
    \right.
    \label{eq:nablan_kappa}
\end{equation}
where we have introduced the principal curvatures $\kappa_i, i=1,\dots, n-1$,
and the Gauss curvature $K = \prod_i \kappa_i$.
The second normal derivative can then be written as:
\begin{equation*}
    \begin{aligned}
        \nabla_{{N}}^2 \kappa 
        & = \nabla_i d \: \nabla_i \nabla_{{N}}\kappa \\
        & = -2 \: \nabla_i d \: \nabla_i \nabla_j \nabla_k d \: \nabla_k \nabla_j d \\
        & = -2 \: \left[\nabla_j \left( \cancel{\nabla_i d \: \nabla_k\nabla_i d} \:
        \nabla_j\nabla_k d \right) - \nabla_j\nabla_i d \: \nabla_k\nabla_i d \:
        \nabla_j\nabla_k d - \cancel{\nabla_i d \: \nabla_k\nabla_i d} \:
        \nabla_j\nabla_j\nabla_k d\right] \\
        & = 2 \: \mathrm{Tr}\left[ (\nabla^2 d)^3 \right]
    \end{aligned}
\end{equation*}
that is:
\begin{equation}
    \nabla_{{N}}^2 \kappa = 2\sum_i \kappa_i^3 
    = \left\{
    \begin{aligned}
        & 2\kappa^3
        \\
        & 2\kappa \left(\kappa^2 - 3K\right)
    \end{aligned}
    \right.
    = \left\{
    \begin{aligned}
        & 2 H^3 &&\textrm{ in 2D}
        \\
        & 4 H \left(4 H^2 - 3K\right) &&\textrm{ in 3D}.
    \end{aligned}
    \right.
    \label{eq:nablan2_kappa}
\end{equation}
Inserting \cref{eq:nablan_kappa} and \cref{eq:nablan2_kappa} into
\cref{eq:delta2_d_eq1}, we finally obtain:
\begin{equation}
    \Delta^2 d \equiv 2w 
    = \left\{ 
    \begin{aligned}
        &\Delta_{\Gamma} \kappa + \kappa^3
        \\
        &\Delta_{\Gamma} \kappa + \kappa (\kappa^2 - 4 K) 
    \end{aligned}
    \right.
    = \left\{ 
    \begin{aligned}
        &\Delta_{\Gamma} H + H^3 &&\textrm{ in 2D}
        \\
        &2 \left( \Delta_{\Gamma} H + 4 H (H^2 - K) \right) &&\textrm{ in 3D}.
    \end{aligned}
    \right.
    \label{eq:Delta2_d}
\end{equation}

As expected, we observe that the second-order term in the expansion of the heat
equation solution \labelcref{eq:Gt_expansion_2} \review{features the required} high-order
geometrical quantities for the geometrical flows we are interested in. In the
following, we shall in particular devise specific algorithms for the simulation
of Willmore flows, but the method presented here is generic, and can be carried
on at higher orders similarly.

\subsection{Convolution of a generic function of the distance}

For practical reasons that will prove useful thereafter, we consider in this
section the expansion of the heat equation solution \cref{eq:heat_equation}, but
with the initial condition:
\begin{equation*}
    \phi(t=0) = f(d)
\end{equation*}
where $f$ is some generic $\mathcal{C}^{\infty}\left(\mathbb{R}\right)$ function.
Similarly to \cref{eq:Gt_expansion_2}, we denote $\mathcal{G}_{\delta t}$ the
corresponding solution, and write as before:
\begin{equation}
    \mathcal{G}_{\delta t}^{(f)} = f(d) + \Delta f(d) \: \delta t + \Delta^2 f(d) \: \frac{\delta
t^2}{2} + \mathcal{O}\left(\delta t^3\right).
    \label{eq:Gtf_expansion1}
\end{equation}

Recalling that $|\nabla d| = 1$, and using the usual derivation rules, we can
then compute:
\begin{equation*}
    \Delta f(d) = \nabla \cdot (f^{(1)} \nabla d) = f^{(2)} + f^{(1)} \Delta d
\end{equation*}
and:
\begin{equation*}
    \begin{aligned}
        \Delta^2 f(d)
        & = f^{(4)}(d) + 2f^{(3)}(d) \Delta d + f^{(2)}(d) \left[ (\Delta d)^2 +
        2\nabla_{{N}} \Delta d \right] + f^{(1)}(d) \Delta^2 d
        \\
        & = f^{(4)}(d) + 2f^{(3)}(d) \kappa + f^{(2)}(d) \left[ \kappa^2 +
        2 \nabla_{{N}} \kappa \right] + 2 f^{(1)}(d) w
        \\
        & \equiv f^{(4)}(d) + 2f^{(3)}(d) \kappa + 2 f^{(2)}(d) c + 2 f^{(1)}(d) w
    \end{aligned}
\end{equation*}
where we have introduced:
\begin{equation}
    c \equiv \frac{\kappa^2}{2} + \nabla_{{N}} \kappa 
    = \left\{ 
    \begin{aligned}
        & -\frac{\kappa^2}{2}
        \\
        & 2K - \frac{\kappa^2}{2}
    \end{aligned}
    \right.
    = \left\{ 
    \begin{aligned}
        & -\frac{H^2}{2} &&\textrm{ in 2D}
        \\
        & -2 (H^2 - K) &&\textrm{ in 3D}.
    \end{aligned}
    \right.
    \label{eq:correction_term}
\end{equation}
We eventually get:
\begin{equation}
    \mathcal{G}_{\delta t}^{(f)} = 
    f(d) 
    + \delta t  \left( f^{(2)}(d) + f^{(1)}(d) \kappa \right)
    + \delta t^2 \left( \frac{f^{(4)}(d)}{2} + f^{(3)}(d) \kappa + f^{(2)}(d) c + f^{(1)}(d) w \right)
    + \mathcal{O}\left(\delta t^3\right).
    \label{eq:Gfdt}
\end{equation}

\subsection{Expansion order selection}
\label{subsec:ExpansionOrderSelection}

We now turn to the use of combinations of $\mathcal{G}$ for different time-steps
in order to select the appropriate order in the expansion
\cref{eq:Gt_expansion_2}. To minimise assembly and
preconditioner computation costs, we restrict our analysis to the case of two
different time-steps denoted $\alpha \delta t$ and $\beta \delta t$
respectively. To simplify the notations, we denote in the following 
$\mathcal{G}_{\alpha} \equiv \mathcal{G}_{\alpha\delta t}$
and
$\mathcal{G}_{\beta} \equiv \mathcal{G}_{\beta\delta t}$.
We then consider the general linear combination 
$a \, \mathcal{G}_{\alpha} + b \, \mathcal{G}_{\beta}$
\review{which upon expansion reads}:
\begin{equation*}
    a \, \mathcal{G}_{\alpha} + b \, \mathcal{G}_{\beta} = 
    (a+b)\: d + (a\alpha+b\beta)\: \kappa \delta t + (a\alpha^2+b\beta^2)\: w \delta
    t^2
\end{equation*}

\reviewII{The coefficients can then be
adjusted to cancel the term of order $\delta t^2$ (the ``Order~1''
combination) or the term of order $\delta t$ (the ``Order~2'' combination).
As mentioned above, we want to use only $\alpha$ and $\beta$ for both
combinations; but the multiplicative coefficients $a$ and $b$ can be chosen
arbitrarily. We thus introduce the ``Order~1'' (resp. ``Order~2'') coefficients
$a_1$ and $b_1$ (resp. $a_2$ and $b_2$) and write the corresponding order
selections as:}
\begin{alignat}{3}
    &\mathrm{Order \ 1: }\quad & a_1 \alpha+b_1 \beta &= 1, & \quad a_1 \alpha^2+b_1 \beta^2 &= 0 \\
    &\mathrm{Order \ 2: }\quad & a_2 \alpha+b_2 \beta &= 0, & \quad a_2 \alpha^2+b_2 \beta^2 &= -1
\end{alignat}
\reviewII{Note that we have arbitrarily imposed $a_1 \alpha+b_1 \beta = 1$ and
    $a_2 \alpha^2+b_2 \beta^2 = -1$ (these quantities actually only need to be non-zero)
    to directly recover the combinations used for our Willmore flow algorithms
below.}

With the additional constraint that $\alpha > 0$, $\beta > 0$, the general
solution is:
\begin{equation*}
    \beta = \alpha \frac{\sqrt{-a_1}}{\sqrt{b_1}}, \quad b_2 = a_2
    \frac{\sqrt{b_1}}{\sqrt{-a_1}}. 
\end{equation*}
We thus choose the following parameters:
\begin{equation*}
    \begin{aligned}
        & \alpha = \sqrt{2}, \quad \beta = \frac{\sqrt{2}}{2}
        \\
        & a_1 = -\frac{\sqrt{2}}{2}, \quad b_1 = 2\sqrt{2}
        \\
        & a_2 = -1, \quad b_2 = 2
    \end{aligned}
\end{equation*}
so that:
\begin{equation}
    a_1 \, \mathcal{G}_{\alpha} + b_1 \, \mathcal{G}_{\beta} =
    \frac{3\sqrt{2}}{2} d + \kappa \delta t + \mathcal{O}\left(\delta t^3\right)
    \label{eq:order1_combination}
\end{equation}
and:
\begin{equation}
    a_2 \, \mathcal{G}_{\alpha} + b_2 \, \mathcal{G}_{\beta} =
    d - w \delta t^2 + \mathcal{O}\left(\delta t^3\right)
    \label{eq:order2_combination}
\end{equation}

\subsection{Willmore flow}
\label{Sec:WillmoreFlow}
In this section, we derive an algorithm which relies on solving heat equation
problems starting from a distance function to compute high-order motions of an
interface in a two- or three-dimensional domain. We focus our analysis more
specifically to the Willmore flow, but the approach can be straightforwardly
extended to a large class of surface diffusion motions.

We consider the Willmore energy functional \cite{willmore1996riemannian}:
\begin{equation}
    \mathcal{E}_W \left(\Gamma\right) \equiv \int_{\Gamma} H^2.
    \label{eq:willmore_energy}
\end{equation}
The surfaces which minimise \cref{eq:willmore_energy} over the set of one- or two-dimensional
surfaces can be obtained from the \textit{Willmore flow} gradient descent, which
moves iteratively
some initial surface with the normal velocity $W {N}$ where (c.f.
\cite{willmore1996riemannian}):
\begin{equation*}
    W = \left\{
        \begin{aligned}
            &\Delta_{\Gamma} H + \frac{H^3}{2} &&\textrm{ in 2D}
            \\
            &\Delta_{\Gamma} H + 2H (H^2 - K)  &&\textrm{ in 3D}.
        \end{aligned}
    \right.
\end{equation*}

For a surface $\Gamma$ implicitly defined by the signed distance function as
$\Gamma = \left\{ \vec{x}\in \mathbb{R}^3, d(\vec{x}) = 0 \right\}$, the
Willmore flow corresponds to a usual transport equation:
\begin{equation}
    \partial_t d + W N \cdot \nabla d = 0.
\end{equation}
Recalling that $\nabla d = N$, this amounts to solving 
$\partial_t d + W =0$, or if we consider a first-order time-discretisation with
time-step $\delta t$,
\begin{equation}
    d^{(p+1)} = d^{(p)} - W \delta t.
    \label{eq:d_transport_explicit}
\end{equation}

We can then observe that the right-hand side of \cref{eq:d_transport_explicit}
can be obtained from the expansions of solutions to the heat equation derived
above. More precisely, the second order term $w\delta t^2$ (\cref{eq:Delta2_d})
provides the required high-order surface term $\Delta_{\Gamma} H$. This is
reminiscent of the approach of 
\cite{esedoglu2008threshold,grzhibovskis2008convolution}, where convolutions of
the characteristic function defining the interface are combined to evolve the
interface according to the Willmore flow. However, while the second-order terms
appearing in the expansion of the convolution of the characteristic function by
the heat kernel provide the exact Willmore velocity, this is no more the case
for the distance function, and we need to add some correction terms in order to
retrieve the correct Willmore flow.

We also stress that $d^{(n)} - W \delta t$ is not in general a distance
function, so that we need to add a redistanciation step to iterate the process.
This redistanciation is the counterpart to the thresholding step in the usual
convolution-thresholding algorithms for mean curvature
\cite{merriman1992diffusion} or higher-order
\cite{esedoglu2008threshold,grzhibovskis2008convolution}
motions using a characteristic function. In the following, we denote
$\mathrm{\bf Redist}$ this redistanciation step, and will use a fast-marching
algorithm to actually compute it. We provide more details about our
implementation in \cref{sec:Numerical_section}.

\paragraph{2D case} In 2D, the correction is rather simple, and only involves
the curvature $H^3$, namely:
\begin{equation}
    W_{2D} = 2w - \frac{H^3}{2} = \frac{2}{\delta t^2} \left( d - a_2 \,
    \mathcal{G}_{\alpha} - b_2 \, \mathcal{G}_{\beta} \right) -
    \frac{H^3}{2}.
    \label{eq:W2d}
\end{equation}

The curvature $H$ can be computed with different methods, including using the
order-$1$ selection \cref{eq:order1_combination}, the $c$ correction term
introduced above in \cref{eq:correction_term} or direct methods using $\Delta
d$. In our computations we will use either the order-$1$ diffusion combination
method or the $c$ correction term computed with diffusions of $d^2$ (c.f.
\cref{eq:Gd2_order2}) as for the
three-dimensional case. 
We refer to \cref{subsec:curvature_computation} for a discussion on the efficiency 
of each method.

As pointed out in \cref{rk:rem_time_discret} and \review{discussed in more detail} in
\cref{solving_the_diff_eq}, the choice of the numerical scheme used to solve the heat
equation has an impact on the development \cref{eq:Gt_expansion_1}. We present
here algorithms using either an implicit Euler or a Crank-Nicolson scheme and choose 
our combination accordingly.

Our diffusion-redistanciation algorithm for the Willmore flow of a curve in two
dimensions is finally given in \cref{alg:2dWillmoreFlow}.

\begin{algorithm}
    \caption{2D Willmore flow\cref{note:algTypo}}
    \label{alg:2dWillmoreFlow}
    \begin{algorithmic}
        \WHILE{$t < t_f$}
        \STATE{Solve $\partial_t \phi - \Delta \phi = 0$ with initial condition
            $\phi_i = d^{(n)}$ for times $t=\sqrt{2}\delta t$ and
            $t=\delta t/\sqrt{2}$
            $\rightarrow \mathcal{G}_{\sqrt{2}\delta t}, \: \mathcal{G}_{\delta t/\sqrt{2}}$. }
            \STATE{Compute the curvature $H$.}
            \STATE{Compute $\mathcal{D}=
			\left\lbrace
			\begin{aligned}            
            & \mathcal{D}_{E} = - \mathcal{G}_{\sqrt{2}\delta t} +
            2\mathcal{G}_{\delta t/\sqrt{2}} + \frac{H^3}{2} \delta t^2 \\
            & \mathcal{D}_{\textsc{CN}} = - \mathcal{G}_{\sqrt{2}\delta t} +
            2\mathcal{G}_{\delta t/\sqrt{2}} + \frac{H^3}{4} \delta t^2 \end{aligned} \right.$
            .}
            \STATE{Compute the new signed distance function $d^{(n+1)} =
            \mathrm{\bf Redist}(\mathcal{D})$.}
            \STATE{$t \leftarrow \left\lbrace
            \begin{aligned} 
                &t + \delta t^2 &\quad\textrm{(Euler)}\\
                &t + \frac{\delta t^2}{2} &\quad\textrm{(Crank-Nicolson)}
            \end{aligned}\right.$
            }
        \ENDWHILE
    \end{algorithmic}
\end{algorithm}
During our simulations, we however preferred to use an algorithm that computes
the correction the same as in the 3D case method using the convolution of $d^2$,
presented in the next paragraph, as it proved more robust in some cases. 
The corresponding algorithm is described in \cref{alg:2dWillmoreFlowWithD2}.
\begin{algorithm}
    \caption{2D Willmore flow with alternative correction\cref{note:algTypo}}
    \label{alg:2dWillmoreFlowWithD2}
    \begin{algorithmic}
        \WHILE{$t < t_f$}
        \STATE{Solve $\partial_t \phi - \Delta \phi = 0$ with initial condition
            $\phi_i = d^{(n)}$ for times $t=\sqrt{2}\delta t$ and
            $t=\delta t/\sqrt{2}$
            $\rightarrow \mathcal{G}_{\sqrt{2}\delta t}, \: \mathcal{G}_{\delta t/\sqrt{2}}$. }
            \STATE{Solve $\partial_t \phi - \Delta \phi = 0$ with initial condition
            $\phi_i = d^{(n)2}$ for times $t=\sqrt{2}\delta t$ and $t=\delta t/\sqrt{2}$
            $\rightarrow \mathcal{G}^{d^2}_{\sqrt{2}\delta t}, \: \mathcal{G}^{d^2}_{\delta t/\sqrt{2}}$.}		
            \STATE{Compute the curvature $H$.}
            \STATE{
                Compute $\mathcal{D} = \frac{H d^{(n)2}}{4} 
                + \left(-\mathcal{G}_{\sqrt{2}\delta t} + 2\mathcal{G}_{\delta t/\sqrt{2}}\right) \cdot \left( 1 - \frac{Hd^{(n)}}{2} \right) 
                + \left(-\mathcal{G}^{d^2}_{\sqrt{2}\delta t} + 2\mathcal{G}^{d^2}_{\delta t/\sqrt{2}}\right) \cdot \frac{H}{4}
                $
            }
            \STATE{Compute the new signed distance function $d^{(n+1)} =
            \mathrm{\bf Redist}(\mathcal{D})$.}
            \STATE{$t \leftarrow \left\lbrace
            \begin{aligned} 
                &t + \delta t^2 &\quad\textrm{(Euler)}\\
                &t + \frac{\delta t^2}{2} &\quad\textrm{(Crank-Nicolson)}
            \end{aligned}\right.$
            }
        \ENDWHILE
    \end{algorithmic}
\end{algorithm}
Despite two additional convolutions, this algorithm seems to produce better
results for our numerical experiments. It
has the advantage of using a correction term which is already of order two. We suspect
that multiplying the correction by $\delta t^2$ in the first version can have some
detrimental numerical effects due to the space discretisation errors.

\paragraph{3D case} In 3D, the correction term is more complicated, as it
involves both the mean and the Gaussian curvatures. More precisely, we have:
\begin{equation*}
    W_{3D} = w - 2H(H^2-K) = w + Hc
\end{equation*}
where $c$ is the term introduced in \cref{eq:correction_term}. To compute this
correction term $c$, we note from \cref{eq:Gfdt} that we need to choose $f$ such
that its second derivative does not vanish. To this end, we choose arbitrarily $f(d) = d^2$ and introduce
\begin{equation}
    \mathcal{G}_{\alpha}^{d^2} = d^2
    + 2 \alpha\delta t  \left( 1 + d \kappa \right)
    + 2 (\alpha\delta t)^2 \left( c + d w \right)
    + \mathcal{O}\left(\delta t^3\right).
    \label{eq:Galpha_d2}
\end{equation}
so that we can extract the correction term from $H \mathcal{G}_{\alpha}^{d^2}$.

However, using only $\mathcal{G}_{\alpha}^{d^2}$ to compute $c$ requires to
compensate the zero- and first-order terms in \cref{eq:Galpha_d2}, which can
lead to large errors in particular regarding the $2\delta t \, d \, \kappa$ term. Note
that while this latter term disappears as $h\to 0$ on the interface (as
$d|_{\Gamma}=0$ by definition), the discrete term brings up
$\mathcal{O}\left(h\right)$ errors which can be large when in the $\delta t$
order term. 
To prevent this, we choose to also select directly the second-order term in
\cref{eq:Galpha_d2} using a linear combination as introduced in \cref{eq:order2_combination}. We thus use:
\begin{equation}
    a_2 \, \mathcal{G}_{\alpha}^{d^2} + b_2 \, \mathcal{G}_{\beta}^{d^2}
    = 
    d^2 - 2 \delta t^2 \left( c + dw \right)
    + \mathcal{O}\left(\delta t^3\right).
    \label{eq:Gd2_order2}
\end{equation}
Recalling $W_{3D} = w + Hc$, and combining \cref{eq:order2_combination} and
\cref{eq:Gd2_order2}, we finally obtain:
\begin{equation}
    W_{3D} = \frac{1}{\delta t^2} \left[
        d - \frac{Hd^2}{2} - \left( a_2 \, \mathcal{G}_{\alpha} + b_2 \, \mathcal{G}_{\beta}\right) (1-H\,d) 
        - \frac{H}{2} \left( a_2 \, \mathcal{G}_{\alpha}^{d^2} + b_2 \, \mathcal{G}_{\beta}^{d^2} \right)
    \right]
    \label{eq:W3D}
\end{equation}
and our algorithm for the \review{three-dimensional Willmore flow is given in
\cref{alg:3dWillmoreFlow}}.
\begin{algorithm}
    \caption{3D Willmore flow\cref{note:algTypo}}
    \label{alg:3dWillmoreFlow}
    \begin{algorithmic}
        \WHILE{$t < t_f$}
            \STATE{Solve $\partial_t \phi - \Delta \phi = 0$ with initial condition
            $\phi_i = d^{(n)}$ for times $t=\sqrt{2}\delta t$ and $t=\delta t/\sqrt{2}$
            $\rightarrow \mathcal{G}_{\sqrt{2}\delta t}, \: \mathcal{G}_{\delta t/\sqrt{2}}$.}
            \STATE{Solve $\partial_t \phi - \Delta \phi = 0$ with initial condition
            $\phi_i = d^{(n)2}$ for times $t=\sqrt{2}\delta t$ and $t=\delta t/\sqrt{2}$
            $\rightarrow \mathcal{G}^{d^2}_{\sqrt{2}\delta t}, \: \mathcal{G}^{d^2}_{\delta t/\sqrt{2}}$.}
            \STATE{Compute the curvature $H$.}
            \STATE{
                Compute $\mathcal{D} = \frac{H d^{(n)2}}{2} 
                + \left(-\mathcal{G}_{\sqrt{2}\delta t} + 2\mathcal{G}_{\delta t/\sqrt{2}}\right) \cdot \left( 1 - Hd^{(n)} \right) 
                + \left(-\mathcal{G}^{d^2}_{\sqrt{2}\delta t} + 2\mathcal{G}^{d^2}_{\delta t/\sqrt{2}}\right) \cdot \frac{H}{2}
                $
            }
            \STATE{Compute the new signed distance function $d^{(n+1)} =
                \mathrm{\bf Redist}(\mathcal{D})$.}
            \STATE{$t \leftarrow \left\lbrace
            \begin{aligned} 
                &t + 2\delta t^2 &\quad\textrm{(Euler)}\\
                &t + \delta t^2 &\quad\textrm{(Crank-Nicolson)}
            \end{aligned}\right.$
            }
 
        \ENDWHILE
    \end{algorithmic}
\end{algorithm}

\begin{remark}
    \label{rk:timeStepDecoupling}
    Note that the flow algorithms actually compute a gradient descent with step $\delta
    t^2$. However, the evolution and diffusion time-steps can easily be
    decoupled involving the ratio $\frac{\delta t_{flow}}{\delta t_{diff}^2}$,
    but our numerical tests suggest that keeping this ratio close to $1$ leads
    to the more stable simulations. This decoupling can nevertheless prove
    useful to implement efficient backtracking methods for the choice of the
    flow time-step while avoiding unnecessary resolutions of the diffusion
    equation.
\end{remark}

\section{Diffusion-redistancing schemes with volume and area conservation}
\label{areavolcons}
Since we are interested in computing shapes in $\R^n$, $n=2$ or $n=3$, we will
in this chapter call \textit{volume} the measure with respect to the Lebesgue
measure of $\R^n$ and \textit{area} the surface measure of hypersurfaces of $\R^n$.
\subsection{Review of existing methods}
Problems involving area or volume conservation were addressed in the framework
of diffusion-redistanciation schemes by \cite{kublik2011algorithms}, and their
convergence properties in \cite{lauxswartz2016convergence}. In
\cite{kublik2011algorithms}, a new algorithm for area preserving flows in two
dimensions is introduced by considering normal velocities
$$v_N=\kappa-\bar{\kappa}+S$$ where $\kappa$ is the curvature, $\bar{\kappa}$
its average on the interface, and $S$ is an additional term which depends on the
application. The algorithm relies in 2D on a property linking the mean
curvature with the surface area and the genus number, which in case of a connected
hypersurface $S$ of genus $0$ enclosing an open bounded set $\Om$, boils down
to $$\bar{\kappa}=\frac{2\pi}{|S|}.$$ Using the divergence theorem, one can then
compute $|S|$ with a volume integral as $$|S|=\int_\Om\Delta d\,dx.$$ In the
case where $S=0$, to ensure an accurate volume conservation, a Newton method is
applied to find a real number $\lambda^*$ close to $\bar{\kappa}$ computed above
(which is used as initialisation) to correct the diffusion generated motion by
raising or lowering the convolution of the signed distance function (as was
introduced the framework of convolution-thresholding schemes by
\cite{ruuth2003simple}). In the case $n=3$, while the trick to compute
efficiently the mean curvature is not anymore valid, we could still apply the
method by computing numerically the mean of the mean curvature.

The level-set community has also developed several techniques to address volume
conservation \cite{sussman1998improved, sussman1999efficient,russo2000remark} in
the modelling of multiphase flows. This corresponds to correct the
Hamilton-Jacobi equation of the redistanciation step by a term which ensures
that conservation. Namely, the following equation: $$\partial_\tau\phi =
\operatorname{sgn}(\phi_0)(1-\|\nabla\phi\|) + \lambda f(\phi)$$ is solved with
$\lambda$ computed so that $\int_{\Om_{ij}}H(\phi)dx$ is conserved on each cell
$\Om_{ij}$ of the grid, and $f$ chosen to localise around the interface.

An even more straightforward approach is a post- or
preprocessing trick due to Smolianski \cite{smolianski2005finite} in the level-set
framework, where a raising parameter for the level-set function is explicitly
computed to restore the target volume. The trick relies on an expansion of the volume enclosed by a level-set
in terms of the height of this level-set:
$$|\{\phi - \delta c < 0\}|:=\int_{\{\phi<\delta c\}}1 dx = |\{\phi  < 0\}|+\delta c \int_{\{\phi=0\}}1 d\sigma+o(\delta c)$$
Therefore starting from a reference volume $V_0$, one can define  $\delta c$ so that  $|\{\phi - \delta c < 0\}|=V_0+o(\delta c)$ by setting:
$$\delta c=\frac{V_0-|\{\phi  < 0\}|}{|\{\phi = 0\}|}.$$

Regarding the conservation of the surface area, one approach, proposed in the
context of vesicle or red blood cell simulations is to relax this
conservation by introducing an area change energy with a high stiffness. In
\cite{cottet2004level,cottet2006level} the authors showed that such an energy
could be expressed thanks to $\|\nabla \phi\|$ which records the area change of
level-sets of $\phi$ when $\phi$ is advected by a divergence free vector field.
By taking a high areal tension stiffness, $\|\nabla \phi\|$ is kept close to $1$ in a
neighbourhood of the interface. In practice, this however
induces a high stiffness of the numerical method, which leads to severe time step
constraints. Another approach is to enforce exact zero surface divergence of the
velocity field which advects the level-set function. This can be done using
Lagrange multipliers \cite{doyeux2013simulation, laadhari2014computing}, leading
to good area conservation, but at the expense of a bad conditioning of the
underlying linear systems to be solve at each iterations, which increases the
computational cost.

In this work, we introduce a correction in the spirit of Smolianski's trick to
conserve both the volume and the area. From an optimisation point of view, our
method can be regarded as a projection method, and the whole constrained Willmore
flow as some kind of projected gradient method.

\subsection{Raising a level-set function to achieve some given area and enclosed volume}
Obviously, we cannot in general \review{fulfill} both the area and volume constraints by
adding a constant to a level-set function. So let us look at the case where we
would add a non-constant function $\delta c:\Om\to\R$ to an arbitrary level-set
function $\phi$.
Let $\Hv$ denote a one dimensional smoothed Heaviside function and $\zeta$ its derivative. A typical choice is: 
\begin{equation}
    \Hv(r)=
    \begin{cases}
        \frac12\left(r+1+\frac1\pi\sin(\pi r)\right)&|r|<1\\0&r\le -1\\1&r\ge1
    \end{cases}
    \text{ and } 
    \zeta(r)=
    \begin{cases}
        \frac12\left(1+\cos(\pi r)\right)&|r|<1\\0&|r|\ge1
    \end{cases}
\end{equation}

We can then compute the enclosed volume of the raised level-set:
\begin{equation}
    \label{A3:eqn:rescDevelVol}
    \begin{split}
        |\{\phi - \delta c < 0\}| 
        &= \lim_{\varepsilon\to0}\int_\Om \Hv\left( -\frac{\phi-\delta c}{\varepsilon} \right) dx \\
        &= \lim_{\varepsilon\to0}\int_\Om  \Hv\left(-\frac{\phi}{\varepsilon} \right) dx 
        + \int_\Om \frac{1}{\varepsilon} \zeta\left( \frac{\phi}{\varepsilon}
        \right) \delta c\; dx + \mathcal{O}\left( \delta c^2 \right)
    \end{split}
\end{equation}

Note that if one adds to a level-set function a $\delta c$ which
averages to zero on $\{\phi=0\}$, the enclosed volume does not change. 

\review{Considering now the area, we} have to take care that
$\phi-\delta c$ is not in general a distance function:

\begin{equation*}
    \begin{split}
        |\{\phi - \delta c = 0\}| 
        &=  \lim_{\varepsilon\to0}\int_\Om \frac{1}{\varepsilon} \zeta\left( \frac{\phi-\delta c}{\varepsilon}\right) |\nabla(\phi-\delta c)|dx + \mathcal{O}\left( \delta c^2 \right)\\
        &=  \lim_{\varepsilon\to0}\int_\Om \frac{1}{\varepsilon} \zeta\left( \frac{\phi}{\varepsilon}\right) |\nabla\phi| dx - \int_\Om \frac{1}{\varepsilon^2} \zeta'\left( \frac{\phi}{\varepsilon}\right) |\nabla\phi|  \delta c\; dx - \int_\Om \frac{1}{\varepsilon}\zeta\left(\frac{\phi}{\varepsilon} \right)\frac{\nabla\phi\cdot\nabla \delta c}{|\nabla\phi|} dx + \mathcal{O}\left( \delta c^2 \right) \\
    \end{split}
\end{equation*}
Integrating by parts in the last term:
\begin{equation*}
    \begin{split}
        - \int_\Om \frac{1}{\varepsilon}\zeta\left( \frac{\phi}{\varepsilon}\right) \frac{\nabla\phi \cdot\nabla \delta c}{|\nabla\phi|} dx & = \int_\Om \div\left( \frac{1}{\varepsilon} \zeta\left( \frac{\phi}{\varepsilon} \right) \frac{\nabla\phi}{|\nabla\phi|} \right)\delta c\; dx \\
        & = \int_\Om \frac{1}{\varepsilon^2}\zeta'\left( \frac{\phi}{\varepsilon} \right) \frac{\nabla \phi \cdot \nabla \phi}{|\nabla\phi|} \delta c\;  dx + \int_\Om \frac{1}{\varepsilon}\zeta\left( \frac{\phi}{\varepsilon} \right) \div\left( \frac{\nabla\phi}{|\nabla\phi|} \right) \delta c\;  dx\\
        & = \int_\Om \frac{1}{\varepsilon^2} \zeta'\left(
        \frac{\phi}{\varepsilon}\right) |\nabla\phi|  \delta c\;  dx + \int_\Om
        \frac{1}{\varepsilon}\zeta\left( \frac{\phi}{\varepsilon} \right)
        \kappa \delta c\;  dx
    \end{split}
\end{equation*}
where $\kappa$ is the total curvature as introduced above. We finally get:
\begin{equation}
    \label{A3:eqn:rescDevelPer}
    |\{\phi - \delta c = 0\}|  
    = |\{\phi  = 0\}|  +  \lim_{\varepsilon\to0} \int_\Om
    \frac{1}{\varepsilon} \zeta\left( \frac{\phi}{\varepsilon}
    \right) \kappa \delta c\; dx + \mathcal{O}\left( \delta c^2\right)
\end{equation}

Note that when $\phi$ is a distance function $\phi = d$, the volume and area of
the perturbed level-set $d - \delta c$ simplify to formulae involving only
``surface'' integrals :
\begin{equation*}
    \begin{split}
        |\{d - \delta c < 0\}| 
        &= |\{d < 0\}| + \int_{\{d=0\}} \delta c\;d\sigma + \mathcal{O}\left( \delta c^2 \right)
    \end{split}
\end{equation*}
which boils down to the Smolianski formula when $\delta c$ is constant, and
\begin{equation*}
    |\{d - \delta c = 0\}|  
    = |\{d = 0\}| + \int_{\{d=0\}} \kappa \delta c\; d\sigma + \mathcal{O}\left(
    \delta c^2 \right).
\end{equation*}

Given a target volume $V_0$ and area $A_0$, we propose the following fast method
to correct a level-set function so that its zero level-set achieves those targets.
\begin{enumerate}
    \item Raise $\phi$ by a constant $-\lambda$ to achieve the right enclosed volume:  
        \begin{equation}
            \lambda 
            = \frac{V_0-|\{\phi < 0\}| }
            {\int_\Om \frac{1}{\varepsilon} \zeta\left( \frac{\phi}{\varepsilon}\right) dx}
        \end{equation}
    \item Raise $\tilde{\phi} = \phi - \lambda$ by $-\mu(\kappa-\overline{\kappa})$, which does not change
        the enclosed volume, where the constant $\mu$ and $\overline{\kappa}$
        are given by:
        \begin{equation}
            \begin{aligned}
                \mu 
                &= \frac{A_0-|\{\tilde{\phi} = 0\}|}{\int_\Om
                    \frac{1}{\varepsilon} \zeta\left( \frac{\tilde{\phi}}{\varepsilon}\right) 
                \kappa (\kappa-\overline{\kappa})\; dx}
                = \frac{A_0-|\{\tilde{\phi} = 0\}|}{(\overline{\kappa^2}-\overline{\kappa}^2)
                \int_\Om \frac{1}{\varepsilon} \zeta\left( \frac{\tilde{\phi}}{\varepsilon}\right) dx}
                \\
                \overline{\kappa} 
                &= \frac{\int_\Om\frac{1}{\varepsilon} \zeta\left(
                    \frac{\tilde{\phi}}{\varepsilon}\right) 
                \kappa}{\int_\Om\frac{1}{\varepsilon} \zeta\left(
                \frac{\tilde{\phi}}{\varepsilon}\right)}
            \end{aligned}
    \end{equation}
\end{enumerate}

\begin{remark}
    In the practical numerical computations, the ``thickness'' $\varepsilon$ of
    the interface is kept finite and proportional and close to the mesh size
    $h$: $\varepsilon \in [h; 3h]$, following the usual level-set approach.
\end{remark}

Clearly, the method fails if $\overline{\kappa^2}=\overline{\kappa}^2$, which is
the case for a circle (resp. sphere). But in that case the perimeter (resp.
area) and enclosed area (resp. volume) are linked and one cannot set one
independently of the other.

Let us investigate how the former restoring of volume and area change the energy
we are minimising. Assuming an energy of the form:
$$\mathcal{E}[\phi]=\int_{\Om} Ed\nu_\eps$$
where $d\nu_\eps=\ze dx$ and, with a gradient which can be written as $$d\mathcal{E}[\phi](\psi)=\int_{\Om} F \psi d\nu_\eps.$$
\begin{proposition}
    \label{Prop:decreasingEnergy}
    Let us consider a gradient-corrected method to minimise $\mathcal{E}$ iteratively by alternating steepest descent with corrections stages $1$ and $2$ described above. Then, at first order in the descent parameter, the energy does not increase during this iteration.
\end{proposition}
\begin{proof}
    A steepest descent method to minimise this energy amounts to change $\phi$
    to $\phi+\psi$ with $\psi=-\rho F$, which decreases the energy, at first
    order, by $\rho\int_{\Om} F^2 d\nu_\eps$.

    But this motion changed the enclosed volume and the area of interface. The
    volume is now  $\left|\{\phi+\psi<0\}\right|$ which correspond at first
    order to a change of $\rho\int_{\Om}Fd\nu_\eps$. Therefore the $\lambda$
    parameter computed above is equal to:
    $$\lambda = \rho
    \frac{\int_{\Om}Fd\nu_\eps}{\int_{\Om}d\nu_\eps}=:\rho\dashint_\Om F
    d\nu_\eps$$
    where from now on, we will denote with a dashed integral sign the mean value,
    i.e. the integral divided by the measure of the set (with respect to
    $\nu_\eps$) onto which the integral is taken: 
    $$\dashint_{\Om}Fd\nu_\eps=\frac1{\nu_\eps(\Om)}\int_\Om Fd\nu_\eps$$
    The descent motion and this volume correction change the area by:
    $$\rho\int_{\Om}F \kappa d\nu_\eps-\lambda\int_{\Om}\kappa
    d\nu_\eps=\rho\int_{\Om}F (\kappa-\overline{\kappa}) d\nu_\eps.$$ Thus the
    second step, area correction, computes a correction $\nu
    (\kappa-\overline{\kappa}) $ with: $$\mu = \rho \frac{ \dashint_{\Om} F
    (\kappa-\overline{\kappa}) d\nu_\eps}{
    \overline{\kappa^2}-\overline{\kappa}^2}.$$ Those two corrections raise the
    energy by: $$ \lambda \int_{\Om}Fd\nu_\eps+\mu \int_{\Om}F (H-\overline{H})
    d\nu_\eps = \rho\nu_\eps(\Om)\left[ \left(\dashint_{\Om}Fd\nu_\eps\right)^2+
    \frac{ \left(\dashint_{\Om} F (H-\overline{H}) d\nu_\eps\right)^2}{
    \overline{H^2}-\overline{H}^2}\right].$$ 
    Finally, the two successive corrections
    would not raise the energy more than it has been decreased by the descent method
    provided that the following inequality holds:
    $$\left(\dashint_{\Om}Fd\nu_\eps\right)^2+ \frac{ \left(\dashint_{\Om} F
    (H-\overline{H}) d\nu_\eps\right)^2}{ \overline{H^2}-\overline{H}^2}\le
    \dashint_{\Om} F^2 d\nu_\eps.$$ This is the object of the following lemma. We
    therefore justified that our descent-correction method does not increase energy\review{.} 

\end{proof}
\begin{lemma}
    Let $(X,\Sigma,\mu)$ be a measure space, and $A$ a measurable subset of $X$
    with $\mu(A)<+\infty$. Let $f,g\in L^2(A)$ with $\dashint_A gd\mu=0$. Then:
    \begin{equation}
        \left(\dashint_A f d\mu\right)^2\dashint_A g^2
        d\mu+\left(\dashint_A fg d\mu\right)^2\le \dashint_A f^2 d\mu\dashint_A g^2
        d\mu\label{CSM}
    \end{equation}
\end{lemma}
\begin{proof}
    This is an easy extension of the Cauchy-Schwarz inequality. Indeed applying
    the latter for $f-\dashint_A f d\mu$ and $g$ we have: $$\left(\dashint_A
    \left(f-\dashint_A f dx\right)gd\mu\right)^2\le  \dashint_A
    \left(f-\dashint_A f d\mu\right)^2 d\mu\dashint_A g^2 d\mu$$ Since $g$ is of
    zero mean on $A$, we have for the left hand side of the former inequality,
    $$\dashint_A \left(f-\dashint_A f d\mu\right)gd\mu=\dashint_A fgd\mu$$ while
    concerning the right hand side we observe that: $$ \dashint_A
    \left(f-\dashint_A f d\mu\right)^2 d\mu+ \left(\dashint_A f d\mu\right)^2
    =\dashint_A f^2 d\mu$$ which leads to \cref{CSM}.
\end{proof}

\begin{remark}
    \label{rem:modifiedRescalingArea}
    In order to prevent undesired effects from outside of the neighbourhood of
    the \review{interface} we can slightly alter the method during the step recovering
    the area. We rescale with the function $\mu(\kappa-\bar{\kappa})
    \frac{1}{\varepsilon}\zeta\left(\frac{\phi}{\varepsilon}\right)$ rather than
    $\mu(\kappa-\bar{\kappa})$. This amounts to choose $\delta c =
    (\kappa-\bar{\kappa})
    \frac{1}{\varepsilon}\zeta\left(\frac{\phi}{\varepsilon}\right) $ in
    \cref{A3:eqn:rescDevelPer} and to redefine $\bar{\kappa}$ (resp. $\bar{H}$
    as: 
    $$ \bar{\kappa} =
    \frac{\int_\Omega\frac{1}{\varepsilon^2}\zeta^2\left(\frac{\phi}{\varepsilon}\right)
    \kappa}{\int_\Omega
    \frac{1}{\varepsilon^2}\zeta^2\left(\frac{\phi}{\varepsilon}\right)} $$
    In our simulations, we fix the value of the rescaling term $\mu(H
    -\bar{H})$ outside of the neighbourhood of the interface in order to
    have a neutral effect.
\end{remark}

\begin{remark}
    An alternative option is to rescale the volume and the area during the same step
    by solving the following system:
    $$
    \begin{pmatrix} \int_\Omega\frac{1}{\varepsilon}\zeta\left(\frac{\phi}{\varepsilon} \right) & \int_\Omega\frac{1}{\varepsilon}\zeta\left(\frac{\phi}{\varepsilon} \right) H \\ \int_\Omega\frac{1}{\varepsilon}\zeta\left(\frac{\phi}{\varepsilon} \right) H & \int_\Omega\frac{1}{\varepsilon}\zeta\left(\frac{\phi}{\varepsilon} \right) H^2 \end{pmatrix} \cdot \begin{pmatrix} \lambda \\ \mu \end{pmatrix} = \begin{pmatrix} V_0 - |\lbrace\phi < 0 \rbrace | \\ A_0 - |\lbrace \phi = 0  \rbrace| \end{pmatrix} 
    $$
    The first equation (resp. the second) corresponds to the volume expansion
    \cref{A3:eqn:rescDevelVol} (the area expansion \cref{A3:eqn:rescDevelPer})
    with $\delta c = \lambda + \mu H$. It is important to note that the constant
    part $\mu \bar{H}$ of the area term is absorbed inside $\lambda$. Both
    versions are equivalent in theory and yield similar numerical results.
\end{remark}

\subsection{Willmore flow with volume and area conservation constraints}
Combining the algorithms \cref{alg:2dWillmoreFlow,alg:3dWillmoreFlow} for the
Willmore flow derived in the previous section and the above method to recover volume and
area, we can obtain an algorithm to compute the Willmore flow in both dimension
2 and 3 with conservation of volume and area.

The principle is to alternate one diffusion step moving according the Willmore
flow and recovering the volume and area constraints before the redistanciation.
Using \cref{Prop:decreasingEnergy} we can ensure that the energy has globally
decreased at the end of one whole step of the method.

The corresponding algorithms in 2D and 3D are described in
\cref{alg:2dWillmoreFlowWithConservations,alg:3dWillmoreFlowWithConservations}.
\begin{algorithm}
    \caption{2D Willmore flow with constant volume and area\cref{note:algTypo}}
    \label{alg:2dWillmoreFlowWithConservations}
    \begin{algorithmic}
        \WHILE{$t < t_f$}
        \STATE{Solve $\partial_t \phi - \Delta \phi = 0$ with initial condition
            $\phi_i = d^{(n)}$ for times $t=\sqrt{2}\delta t$ and $t=\delta
            t/\sqrt{2}$
            $\rightarrow \mathcal{G}_{\sqrt{2}\delta t}, \: \mathcal{G}_{\delta
        t/\sqrt{2}}$.}
        \STATE{Compute the curvature $H$.}
        \STATE{Compute $\mathcal{D}=
            \left\lbrace
            \begin{aligned}            
                & \mathcal{D}_{E} = -\mathcal{G}_{\sqrt{2}\delta t} +
                2\mathcal{G}_{\delta t/\sqrt{2}} + \frac{H^3}{2} \delta t^2 \\
                & \mathcal{D}_{\textsc{CN}} = - \mathcal{G}_{\sqrt{2}\delta t} +
                2\mathcal{G}_{\delta t/\sqrt{2}} + \frac{H^3}{4} \delta t^2 
            \end{aligned} \right.$
        }

        \STATE{Compute $H(\mathcal{D})$ and its mean $\overline{H}(\mathcal{D})$. }
        \STATE{Compute $\lambda$:
            \begin{equation*}
                \lambda = \frac{V_0 - |\lbrace \mathcal{D} < 0 \rbrace|}{|\lbrace \mathcal{D} = 0 \rbrace|}. 
            \end{equation*}}
        \STATE{Compute $\mu$:
            \begin{equation*}
                \mu = \frac{A_0- |\lbrace \mathcal{D} = 0 \rbrace|}{\int_{\Omega}
                H(H-\overline{H})\frac{1}{\varepsilon}\zeta\left(\frac{\mathcal{D}}{\varepsilon}\right)d\sigma}.
        \end{equation*}}
        \STATE{Construct the new signed distance function $d^{(n+1)} =
            \mathrm{\bf Redist}\left(\mathcal{D} - \lambda - \mu
        (H-\overline{H})\right)$.}
        \STATE{$t \leftarrow \left\lbrace
        \begin{aligned} 
            &t + \delta t^2 &\quad\textrm{(Euler)}\\
            &t + \frac{\delta t^2}{2} &\quad\textrm{(Crank-Nicolson)}
        \end{aligned}\right.$
        }
        \ENDWHILE
    \end{algorithmic}
\end{algorithm}

\begin{algorithm}
    \caption{3D Willmore flow with volume and area constraints\cref{note:algTypo}}
    \label{alg:3dWillmoreFlowWithConservations}
    \begin{algorithmic}
        \WHILE{$t < t_f$}
        \STATE{Solve $\partial_t \phi - \Delta \phi = 0$ with initial condition
        $\phi_i = d^{(n)}$ for times $t=\sqrt{2}\delta t$ and $t=\delta t/\sqrt{2}$
        $\rightarrow \mathcal{G}_{\sqrt{2}\delta t}, \: \mathcal{G}_{\delta t/\sqrt{2}}$.}
        \STATE{Solve $\partial_t \phi - \Delta \phi = 0$ with initial condition
        $\phi_i = d^{(n)2}$ for times $t=\sqrt{2}\delta t$ and $t=\delta t/\sqrt{2}$
        $\rightarrow \mathcal{G}^{d^2}_{\sqrt{2}\delta t}, \: \mathcal{G}^{d^2}_{\delta
            t/\sqrt{2}}$.}
        \STATE{Compute the curvature $H$.}
        \STATE{
                Compute $\mathcal{D} = \frac{H d^{(n)2}}{2} 
                + \left(-\mathcal{G}_{\sqrt{2}\delta t} + 2\mathcal{G}_{\delta t/\sqrt{2}}\right) \cdot \left( 1 - Hd^{(n)} \right) 
                + \left(-\mathcal{G}^{d^2}_{\sqrt{2}\delta t} + 2\mathcal{G}^{d^2}_{\delta t/\sqrt{2}}\right) \cdot \frac{H}{2}
                $
            }
        \STATE{Compute $H(\mathcal{D})$ and its mean $\overline{H}(\mathcal{D})$. }
        \STATE{Compute $\lambda$:
        \begin{equation*}
\lambda = \frac{V_0 - |\lbrace \mathcal{D} < 0 \rbrace|}{|\lbrace \mathcal{D} = 0 \rbrace|}. 
\end{equation*}}
        \STATE{Compute $\mu$:
	\begin{equation*}
	\mu = \frac{A_0- |\lbrace \mathcal{D} = 0 \rbrace|}{\int_{\Omega}
        H(H-\overline{H})\frac{1}{\varepsilon}\zeta\left(\frac{\mathcal{D}}{\varepsilon}\right)d\sigma}.
\end{equation*}}
        \STATE{Construct the new signed distance function $d^{(n+1)} =
            \mathrm{\bf Redist}\left(
        {\mathcal{D} - \lambda - \mu (H-\overline{H}})\right)$.}
        \STATE{$t \leftarrow \left\lbrace
        \begin{aligned} 
            &t + 2\delta t^2 &\quad\textrm{(Euler)}\\
            &t + \delta t^2 &\quad\textrm{(Crank-Nicolson)}
        \end{aligned}\right.$
        }
        \ENDWHILE
    \end{algorithmic}
\end{algorithm}

\section{Numerical method and practical implementation}
\label{sec:Numerical_section}
\subsection{Finite Element library}
\label{subsec:numerical_implementation}

The equations and integrals introduced above are solved or evaluated within a
finite-element framework, using the \textsc{Feel++}--finite-element C++ library
\cite{prud2012feel++,prud2006domain}, and in
particular the \textit{LevelSet} framework from the \textsc{Feel++} toolboxes
\cite{metivet2018high}, which features a comprehensive and seamless parallel set
of tools for this kind of surface capturing methods. \reviewII{This framework was preferred to a more classical finite difference approach for the simplicity of coding with high level templates that mimic the variational form of PDE. Moreover, our final goal is to solve fluid-structure interaction problems where vesicles are evolving inside a fluid which flows in a complex geometry like blood vessels.}

More precisely, we use a continuous Galerkin \review{variational} approach, and
discretise the resulting equations in space with Lagrange polynomials.
Introducing $\mathcal{T}_h \equiv \left\{ K_e, 1\leq e \leq N_{elt} \right\}$ a
compatible tessellation of the computational domain $\Omega$, and the
corresponding discrete -- unstructured -- mesh $\Omega_h =
\bigcup_{e=1}^{N_{elt}} K_e$, we define $\mathcal{P}_{h}^{k}\equiv
\mathcal{P}_{h}^{k}\left( \Omega_h\right)$ as the finite-element space on
$\Omega_h$ spanned by Lagrange polynomials
of order $k$.

\subsection{Solving the diffusion equation}
\label{solving_the_diff_eq}
From a numerical perspective, our Willmore flow
\cref{alg:2dWillmoreFlow,alg:3dWillmoreFlow} mainly involve solving the
diffusion equation $\partial_t \phi - \Delta \phi = 0$ with appropriate initial
and boundary conditions. To this end, we use the $\mathcal{P}_{h}^{k}$ finite-element 
space discretisation introduced above, and a second order \review{unconditionally}
stable Crank-Nicolson scheme for the time discretisation. 
The corresponding discrete variational problem then
\review{reads}:
\begin{equation}
    \begin{aligned}
        &\textit{Find $\phi^{(n+1)}\in \mathcal{P}_{h}^{k}$ s.t.
            $\forall \psi\in \mathcal{P}_{h}^{k}$,}
        \\
        &\int_{\Omega_h} \left( 
            \frac{\phi^{(n+1)}}{\delta t} \, \psi 
            + \frac{1}{2} \, \nabla\phi^{(n+1)} \cdot \nabla\psi
        \right)
        =
        \int_{\Omega_h} \left( 
            \frac{\phi^{(n)}}{\delta t} \, \psi
            - \frac{1}{2} \, \nabla\phi^{(n)} \cdot \nabla\psi
        \right)
        + \int_{\partial\Omega_h} \left( {N} \cdot \nabla\phi^{(n)} \right) \, \psi
    \end{aligned}
    \label{eq:diffeq_variationnal}
\end{equation}
where the superscript indices denote the time iterations. Note that the last term in
\cref{eq:diffeq_variationnal} comes from the integration by parts of the
diffusive terms, and is somehow similar to an explicit discrete Neumann boundary
condition preserving the normal gradient of $\phi$, namely
${N}\cdot\nabla\phi^{(n+1)} = {N}\cdot\nabla\phi^{(n)}$ with ${N}$
the exterior normal at the boundary of $\Omega_h$. In practice, the best choice of
boundary conditions for our distance diffusion problem is an open question,
since no natural condition emerges from our analysis, which focuses on the
neighborhood of the interface $\Gamma$, far from the domain boundary.

As pointed out in \cref{rk:rem_time_discret}, we would also like to use in some
situations an Euler scheme, which is low order and diffusive, which may be more
adapted in this context than a dispersive scheme. In that case one should
directly use the discrete expansion, which differs at order
two (since the scheme is order one):
$$\phi^{n+1}-\delta t\Delta \phi^{n+1}=\phi^n$$
indeed gives
\begin{equation}
    \label{euler_expansion}
    \phi^{n+1} = (\operatorname{id}-\delta t\Delta )^{-1}(\phi^n)
    = \phi^n+\delta t \Delta\phi^n+\delta t^2 \Delta^2\phi^n+o(\delta t^2)
\end{equation}
which differs from of (\ref{eq:Gt_expansion_1}) by a factor $2$ in the second
order term. The corresponding discrete variational problem then
simplifies to:
\begin{equation}
    \begin{aligned}
        &\textit{Find $\phi^{(n+1)}\in \mathcal{P}_{h}^{k}$ s.t.
            $\forall \psi\in \mathcal{P}_{h}^{k}$,}
        \\
        &\int_{\Omega_h} \left( 
            \frac{\phi^{(n+1)}}{\delta t} \, \psi 
            + \, \nabla\phi^{(n+1)} \cdot \nabla\psi
        \right)
        =
        \int_{\Omega_h}
            \frac{\phi^{(n)}}{\delta t} \, \psi         
        + \int_{\partial\Omega_h} \left( {N} \cdot \nabla\phi^{(n)} \right) \, \psi
    \end{aligned}
    \label{eq:diffeq_variationnal_euler}
\end{equation}
In our numerical results, both schemes give qualitatively the same equilibrium
shapes. However the Crank-Nicolson scheme seems more accurate, while the Euler
scheme was preferred in dimension 3 where its diffusive behavior brings more
stability and smoothness of the interface.

\subsection{Computation of the curvature}
\label{subsec:curvature_computation}
Our algorithms for the Willmore flow in both 2D \cref{alg:2dWillmoreFlow} and 3D
\cref{alg:3dWillmoreFlow} require the computation of the curvature of the
surface $\Gamma$ in the correction term. We first present the direct methods
currently used in \textsc{Feel++}. We then introduce our diffusion methods to
compute the curvature, which are essentially the classical
convolution/thresholding algorithms of order 1 and 2. Using the direct methods
as references, we discuss the efficiency of our method and focus on the choice
of the time step selection for the diffusion steps.

\subsubsection{Direct methods}
\label{subsubsec:curvature_direct_methods}
From a usual level-set perspective, the curvature of $\Gamma$ can be computed
directly using the divergence of the level-set function $\Delta d$. For
low-order discretisations however, this requires the use of specific strategies,
as standard finite-element derivation decreases the discretisation polynomial
order by 1, which forbids derivations of degrees higher than the polynomial
order. To circumvent this issue, we use a classical Galerkin projection in
$L^2(\Omega_h)$ to maintain the element in $\mathcal{P}_{h}^{1}$. Namely the
computation of the gradient $\nabla d$ is as following: 
\begin{equation*}
    \begin{aligned}
        &\textit{Find $g\in \mathcal{P}_{h}^{1}$ s.t.
            $\forall v\in \mathcal{P}_{h}^{1}$,}
        \\
        &\int_{\Omega_h} g\cdot v
        =
        \int_{\Omega_h} \nabla d \cdot v
    \end{aligned}
\end{equation*}
We can then obtain the curvature with two successive derivation/projection steps
by computing the divergence of the gradient of the distance function:
\begin{equation}
\label{DiscreteCurvatureProjectionL2}
    \begin{aligned}
        & (1) :  \textit{Find $g\in \mathcal{P}_{h}^{1}$ s.t.
            $\forall v\in \mathcal{P}_{h}^{1}$,}
        \\
        & \ \int_{\Omega_h} g\cdot v
        =
        \int_{\Omega_h} \nabla d \cdot v \\
        & (2) : \textit{Find $H\in \mathcal{P}_{h}^{1}$ s.t.
            $\forall w\in \mathcal{P}_{h}^{1}$,}
        \\
        & \ \int_{\Omega_h} H w
        =
        \int_{\Omega_h} \div g \ w \\
    \end{aligned}
\end{equation}
However the $L^2$ projection fills the missing information of the lower order
elements with noisy values, so that one usually resort to a smoothed $L^2$
projection method, where a small diffusion term $-\eta \Delta H$ is added:
\begin{equation}
\label{DiscreteProjectionSmoothed}
    \begin{aligned}
        & (1) :  \textit{Find $g\in \mathcal{P}_{h}^{1}$ s.t.
            $\forall v\in \mathcal{P}_{h}^{1}$,}
        \\
        & \ \int_{\Omega_h} g\cdot v
        =
        \int_{\Omega_h} \nabla d \cdot v \\
        & (2) : \textit{Find $H\in \mathcal{P}_{h}^{1}$ s.t.
            $\forall w\in \mathcal{P}_{h}^{1}$,}
        \\
        & \ \int_{\Omega_h} H w + \eta \int_{\Omega_h} \nabla H \cdot \nabla w - \eta\int_{\partial\Omega_h} \nabla H \cdot N  w 
        =
        \int_{\Omega_h} \div g \ w. \\
    \end{aligned}
\end{equation}
where the smoothing coefficient $\eta$ \review{is typically taken equals to $0.03h$}. We will refer to \cref{DiscreteCurvatureProjectionL2} as the \textit{$L^2$
projection} and to \cref{DiscreteProjectionSmoothed} as the \textit{Smoothed
projection} method.

\subsubsection{Using the diffusion of the signed distance function}
\label{subsubsec:curvature_diffusion_method}
As mentioned above, we can also retrieve the curvature of the surface
represented with the level-set using the diffusion of the signed distance
function.  Recalling the analytical expansion of the diffusion solution
\cref{eq:Gt_expansion_2}, we observe that the curvature of the surface shows up in
the first order term in $\delta t$, and can thus be retrieved as:
\begin{equation}
    H = \frac{\mathcal{G}_{\delta t}-d}{(n-1)\delta t} + \mathcal{O}(\delta t).
    \label{A4:eqn:MeanCurvatureOrder1}
\end{equation}

As explained in \cref{subsec:ExpansionOrderSelection}, at the expense of
an additional resolution of the diffusion equation, we can improve the order of
accuracy and compute the combination \cref{eq:order1_combination} which
eliminates the second order term. We can then compute $H$ at order $\delta t^2$
as:
\begin{equation}
    H = \frac{a_1 \mathcal{G}_\alpha + b_1 \mathcal{G}_\beta - \frac{3\sqrt{2}}{2}d}{(n-1)\delta t} + \mathcal{O}(\delta t^2)
    \label{A4:eqn:MeanCurvatureOrder2}
\end{equation}

In the following we shall refer to the $\mathcal{O}(\delta t)$ method
\cref{A4:eqn:MeanCurvatureOrder1} as the \textit{Order $1$ diffusion} and the
$\mathcal{O}(\delta t^2)$ method \cref{A4:eqn:MeanCurvatureOrder2} as the
\textit{Order $2$ diffusion}. \review{Keep in mind however}  that these approximation orders
refer to the auxiliary diffusion time-step and not to the usual space
discretisation order.

\begin{figure}[!ht]
    \centering
    \input{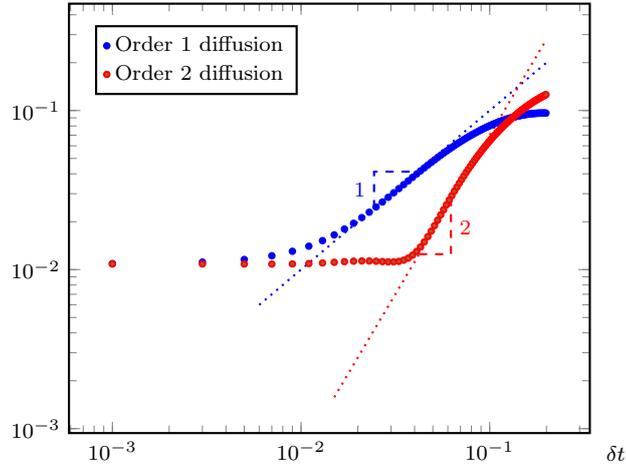}%
    \caption{
        $L^2$ error on curvature estimates using the \textit{Order $1$}
        \cref{A4:eqn:MeanCurvatureOrder1} and \textit{Order $2$}
        \cref{A4:eqn:MeanCurvatureOrder2} diffusion methods as a function of the
        diffusion time-step $\delta t$. We observe the expected convergence
        orders as the time-step decreases, and a saturation of the error for low
        time-steps which is related to space discretisation error terms.
    }
    \label{s4:fig:diffCurvatureError}
\end{figure}

Before comparing these diffusion methods to the direct ones mentioned
\review{in \cref{subsubsec:curvature_direct_methods}}, we shall first study in more
detail the role of the diffusion time-step.  \Cref{s4:fig:diffCurvatureError}
shows the evolution of the $L^2$ error of the curvature computed with both
diffusion methods as a function of the diffusion time-step. The numerical test
was performed on a two-dimensional circular surface of curvature $1$ with a mesh
size $h \approx 0.02$, and the error was computed as 
\begin{equation}
    \mathrm{err}(H;\Gamma) = \frac
    {\int_{\Omega} \left( H - H_{th} \right)^2
    \:\frac{1}{\varepsilon} \zeta\left(\frac{\phi}{\varepsilon}\right)}
    {\int_{\Omega}\frac{1}{\varepsilon}
    \zeta\left(\frac{\phi}{\varepsilon}\right)},
    \label{eq:diffCurvatureL2Err}
\end{equation}
\review{where $H_{th}$ stands for the theoretical value of $H$. In the case of a circle, it amounts to the inverse of the radius.}
\reviewII{The spreading parameter $\varepsilon$ was chosen for all the
    simulations so that the interface
span about $2$ mesh cells, i.e. $\varepsilon \approx h$.}

As expected, we observe convergence of both methods at order $\mathcal{O}(\delta
t)$ and $\mathcal{O}(\delta t^2)$ respectively, as well as a saturation of the
error for small time-steps, which is due to the spatial discretisation errors.
We also note that the \textit{Order $2$ diffusion} method features an optimal
error for larger time-steps, which can prove useful when using the same
time-step for the flow and the diffusions (c.f. \cref{rk:timeStepDecoupling}).

These numerical observations, as well as the analytic expansion of
$\mathcal{G}_{\delta t}$ also suggest to choose the time-step following the
heuristic strategy:
\begin{equation}
    \tilde{H} \delta t \approx h
    \label{OptTimeStepRule}
\end{equation}
with $\tilde{H}$ an \textit{a priori} estimate of the typical curvature of the surface.
From a numerical point of view, this corresponds to having the diffusion move
the $d=0$ level-set a few mesh elements.
This heuristic can also be adapted to the case of geometrical flow simulations,
and naturally formulates an adaptive time-stepping method, as presented in
\cref{subsec:AdaptiveTimeStepping}.

\subsubsection{Comparison of the different methods}
\begin{figure}[!ht]
    \centering
    \input{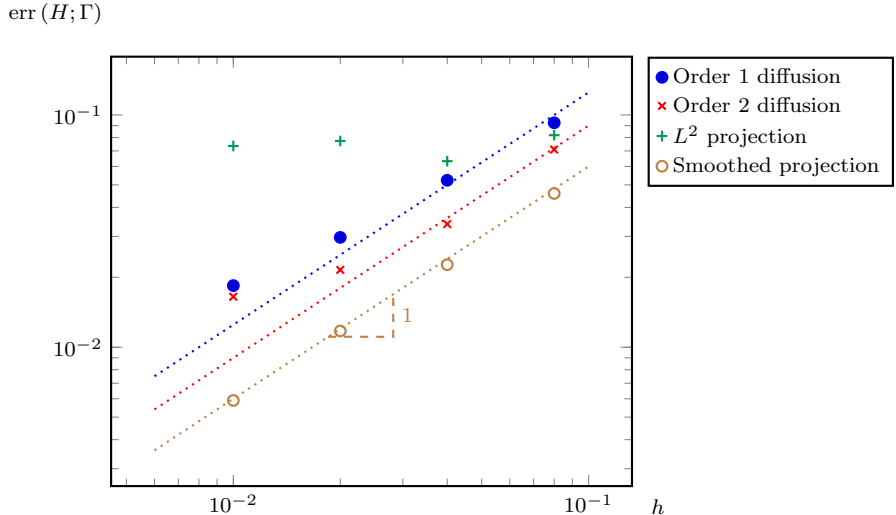}
    \caption{
        Evolution of the $L^2$ error as a function of the mesh size
        for the different methods. The \textit{$L^2$ projection} shows
        no evident convergence while the \textit{Order $1$ diffusion}, the
        \textit{Order $2$ diffusion} and the \textit{Smoothed projection} methods 
        yield linear convergence.
    }
    \label{s4:fig:curvatureMethodsL2Errors}
\end{figure}

We now turn to the comparison of the ``direct'' and ``diffusion'' curvature methods 
presented above.
We again use a fixed two-dimensional circular interface of radius $1$, and vary
the characteristic mesh size $h$ used to compute the $L^2$ error on
the curvature estimate along the surface as a function of $h$ for the four
methods. We set the smoothing parameters $\eta = 0.03 h$ and the diffusion time-steps
as $\delta t = 0.7 h$. The resulting errors are plotted in log-log scale in
\cref{s4:fig:curvatureMethodsL2Errors}.

We find that the \textit{Order $1$} and \textit{Order $2$} diffusion methods behave
similarly, and compare accurately with the \textit{Smoothed projection} method.
These three methods give $\mathcal{O}(h)$ convergence, which is satisfactory for
our low-order discretisation. As anticipated, the \textit{$L^2$ projection} method
does not seem to converge.

In our simulations, we shall therefore use the \textit{Order $2$ diffusion} 
method, as the two diffusion equation resolutions for $\mathcal{G}_\alpha$ and
$\mathcal{G}_\beta$ are needed anyway for the computation of the Willmore flow,
and this method provides the best curvature estimate at no additional cost.

\subsection{Level-set redistanciation}
\label{subsec:levelset-redistanciation}

The redistanciation step is a crucial aspect of our algorithm, as it allows the
iteration to proceed while preserving the correct expansions of the diffusion
solutions close to the interface. This step can moreover become the main
bottleneck of the whole algorithm in terms of both efficiency and stability if
not processed carefully.
From a general point of view, performing the $\mathrm{\bf Redist}$ step in
\cref{alg:2dWillmoreFlow,alg:3dWillmoreFlow} amounts to computing the distance
to an interface $\Gamma$ located by the $0$-level of some function. Stated differently,
we need to solve the eikonal equation:
\begin{equation}
    |\nabla \phi| = 1, \quad \phi(\Gamma) = 0
\end{equation}
which is a challenging boundary-value non-linear equation which has received
much attention since the work of Sethian \cite{sethian1996fast} who proposed an
iterative upwind \textit{fast-marching} algorithm to solve the eikonal equation
starting from the interface and propagating the information outward.

In our finite-element framework, we use a parallel fast-marching
algorithm inspired from \cite{yang2017highly} but adapted to arbitrary --
possibly unstructured -- meshes. The local -- element-wise -- eikonal
equations are solved with a \textit{QR} decomposition, and the fast-marching
alternates between local-domain solves and ghost inter-domain updates until
global convergence is obtained.

In order to ensure good stability properties of our method, and in particular to
prevent the fast-marching method from introducing spurious motion of the interface,
the initialisation of the algorithm near the implicitly defined interface is
crucial. To this end, we rescale the level-set
function $\phi$ by $\frac{1}{|\nabla \phi|}$ on the elements which intersect the
interface, \review{which exactly solves the local eikonal equation when $\phi$
is a piecewise linear polynomial ($\phi\in\mathcal{P}_h^1$)}, before applying the fast-marching algorithm both inward ($\phi < 0$)
and outward ($\phi > 0$) starting from the values encompassing the interface.


\subsection{Adaptive time-stepping}
\label{subsec:AdaptiveTimeStepping}

As can be seen from \cref{eq:order2_combination}, given some $\delta t$ in the
resolution of the heat equations, our diffusion-redistanciation
schemes actually solves the corresponding Willmore flow with time-step $\delta
t^2$, namely
\begin{equation*}
    d^{(n+1)} = d^{(n)} - W \delta t^2.
\end{equation*}
The flow time-step is therefore not constrained by any stability issue, but
controls the accuracy of the solution, since the expansions performed in
\cref{subsec:convolutionexpansion} are valid up to terms of order
$\mathcal{O}(\delta t^3)$. 

The $\mathcal{O}(h)$ errors pertaining to the computation of $\phi^{(n+1)}$ in
\cref{eq:diffeq_variationnal} using $\mathcal{P}_{h}^{1}$ Lagrange elements and
the expansion \labelcref{eq:Gt_expansion_2} however suggest to use a time-step
like
\begin{equation*}
    \delta t \sim \frac{h}{\kappa}
\end{equation*}
in the heat equation. This obviously requires the \textit{a priori} computation of the
interface curvature, which is not desirable in our framework. Instead, we
therefore use an explicit adaptive time-stepping strategy by setting
\begin{equation}
    \delta t^{(n+1)} = \min_{\Gamma}\left(\frac{h}{\kappa^{(n)}}\right)
    \label{s4:eq:adaptiveTimeStep}
\end{equation}
where $\kappa^{(n)}$ is the curvature at the
previous iteration and $h$ is a measure of the mesh size, here considered for
elements crossed by the interface -- e.g.the minimal diameter of the elements crossed by
the $\{\phi^{(n)} = 0\}$ level-set.

\section{Numerical illustrations}
In this section, we illustrate the efficiency our methods by applying them to two classical applications. First, we apply our 3D algorithm for the Willmore flow \cref{alg:3dWillmoreFlow} to the Willmore problem. Then we apply our algorithm for Willmore flow with area and volume conservation \cref{alg:3dWillmoreFlowWithConservations} to the computation of an equilibrium shape of a red blood cell.

\subsection{Willmore problem}
The \textit{Willmore problem} refers to the study of the minimisers of the
Willmore energy function $\mathcal{E}_W$. The case of compact surfaces of genus
$0$ is trivial as the minimisers are the spheres -- recall that the Willmore
energy being scale invariant, all the spheres have the same energy.

In the case of compact surfaces of genus $1$, the Clifford torus, the torus with
a ratio $\sqrt{2}$ between its radii, minimises the Willmore energy. The proof of
this conjecture made by Willmore \cite{willmore1965note} in 1965 has been
established recently in \cite{marques2014min}. The conjectures for higher
genus order \cite{hsu1992minimizing} are still to be proven.

In the following, we present 2D and 3D simulations of such simple surfaces
which evolve according to the Willmore flow using our diffusion-redistanciation
algorithms. The existence of analytic minimisers will then provide solid means
to evaluate our numerical approach.

\subsubsection{2D Willmore flow of a circle}
To assess the accuracy of our diffusion-redistanciation scheme, we first
simulate an initial two-dimensional circular interface evolving according to the
Willmore flow, which can also be computed analytically. As it evolves with the
flow, the interface should remains circular, and increase it radius $r$ as
\begin{equation}
    r(t) = \left( r_0 + 2\, t \right)^{1/4}
    \label{eq:willmoreRadiusSol}
\end{equation}
with $r_0$ the initial radius. Note that this law, which comes from the analytic
flow equation $\frac{dr}{dt} = \frac{H^3}{2} = {1}{2r^3}$ provides a good test
of our algorithm, since the right hand side is actually obtained through a
delicate compensation between the second-order term $w$ (c.f.
\labelcref{eq:Delta2_d}) in the diffusion of $d$ and the correction term $c$
obtained with the diffusion of $d^2$ (c.f. \labelcref{eq:correction_term}).

Our simulation was run with $r_0 = 1$ in a square domain of side length $6$ with 
an unstructured mesh of typical size $h \approx 0.02$ and the
adaptive time-step strategy presented in \cref{subsec:AdaptiveTimeStepping}.

\begin{figure}[!htb]
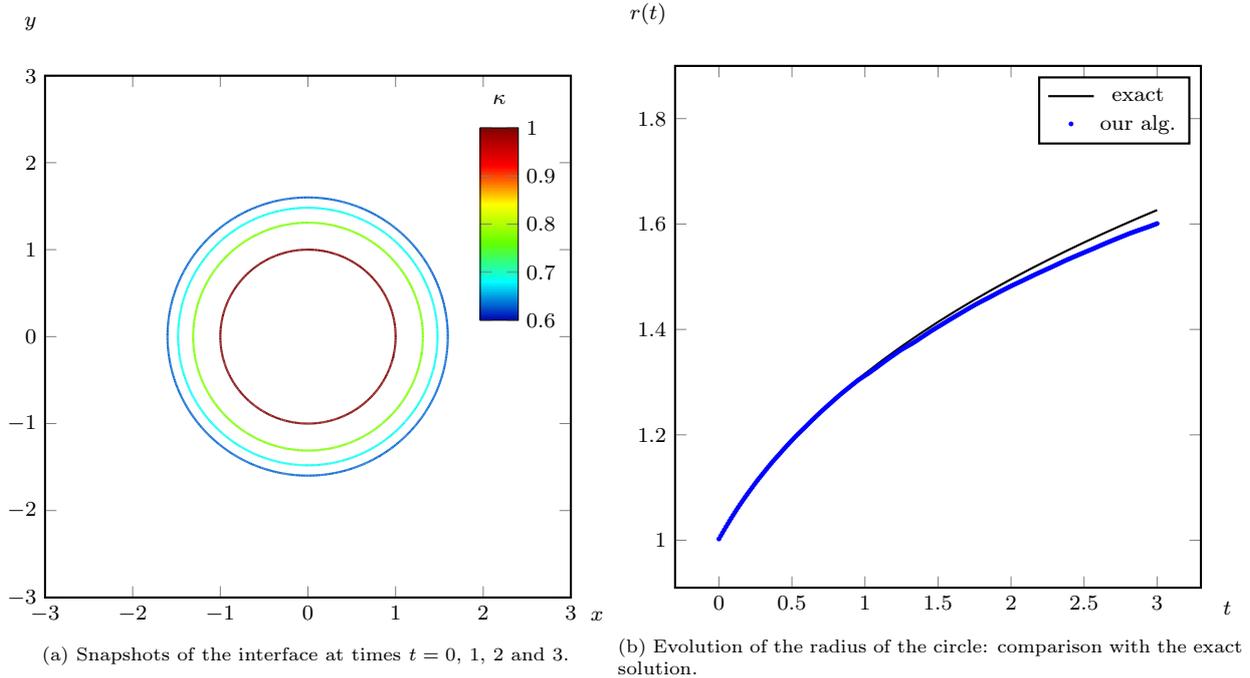

    \centering
    \begin{subfigure}{0.5\textwidth}
        \centering
        \input{figures/willmore_circle2d_evol.tex}%
        \caption{Snapshots of the interface at times $t=0$, $1$, $2$ and $3$.}
    \end{subfigure}%
    \begin{subfigure}{0.5\textwidth}
        \centering
        \input{figures/willmore_circle2d_radius.tex}%
        \caption{Evolution of the radius of the circle: comparison with the
        exact solution.}
    \end{subfigure}%
    \caption{2D Willmore flow of a circle.}
    \label{fig:WillmoreCircle2d}
\end{figure}

\begin{figure}[!ht]
    \centering
    \input{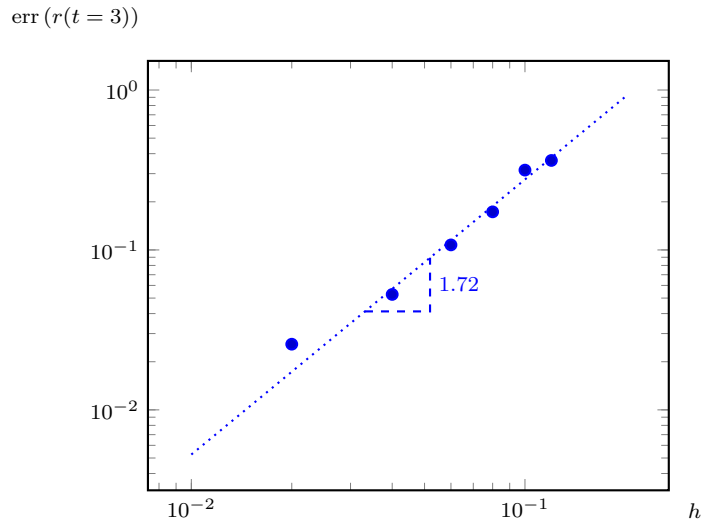}
    \caption{
        Evolution of the error $\mathrm{err}(r(t=3)) \equiv | r_{h}(t=3) - r(t=3) |$
        of the radius of the circle at $t=3$ as a
        function of the mesh size. The figure also shows the best mean-square
        fit of the data by a power function $14.48 \,h^{1.72}$. Note however that
        the displayed slope is not strictly equivalent to the usual convergence
        order, as discussed above.
    }
    \label{s5:fig:willmoreCircle2dRadiusErrors}
\end{figure}

\Cref{fig:WillmoreCircle2d} shows some snapshots of the simulation and the
evolution of the circle radius. It illustrates the numerical accuracy of our
scheme, which succeeds in preserving the symmetry of the shape as it evolves,
and compares quantitatively with the exact solution.
\review{\Cref{s5:fig:willmoreCircle2dRadiusErrors} shows the error on the computed
radius of the circle as compared to the theoretical value
\labelcref{eq:willmoreRadiusSol} at $t=3$ for different mesh sizes. The
dependence of this error on the mesh size demonstrates the convergence of our
method. Note however that the fitted convergence order -- $\approx 1.72$ -- is
not strictly equivalent to the usual finite-element convergence order, as it
incorporates effects related to the dependence of the time step on the mesh
size, which is chosen following \cref{s4:eq:adaptiveTimeStep} as the role played by the 
time step in the derivation of our algorithm naturally imposes such
spatio-temporal coupling. The linear
dependence of the time-step with respect to the mesh size would suggest
that the actual spatial convergence order is $\sim 0.7-0.8$, which also seems
supported by the energy convergence shown in
\cref{s5:fig:vesicle3dEnergiesConv}.}

\subsubsection{3D Willmore flow of a torus}

The convergence of a \review{torus towards a Clifford torus}, whose ratio between major and
minor radii equals $\sqrt{2}$, is a good test for our numerical approach, as it
challenges both the accuracy of the numerical flow and the stability of the
algorithm as it approaches the minimising surface, which features a rather small
hole, which can easily be filled by an inaccurate flow to minimise the energy
further, as the Willmore energy of any sphere -- $4\,\pi$ -- falls below the
minimal one for tori, which is $2\,\pi^2$.

Simulations have been performed using a phase-field model in
\cite{bretin2015phase} but the results are only qualitative.  To the extend of
our knowledge, there is no other published numerical work to compare with our
results.

\begin{figure}[!htb]
    \centering
    \begin{subfigure}{0.33\textwidth}
        \centering
        \includegraphics[width=\textwidth]{{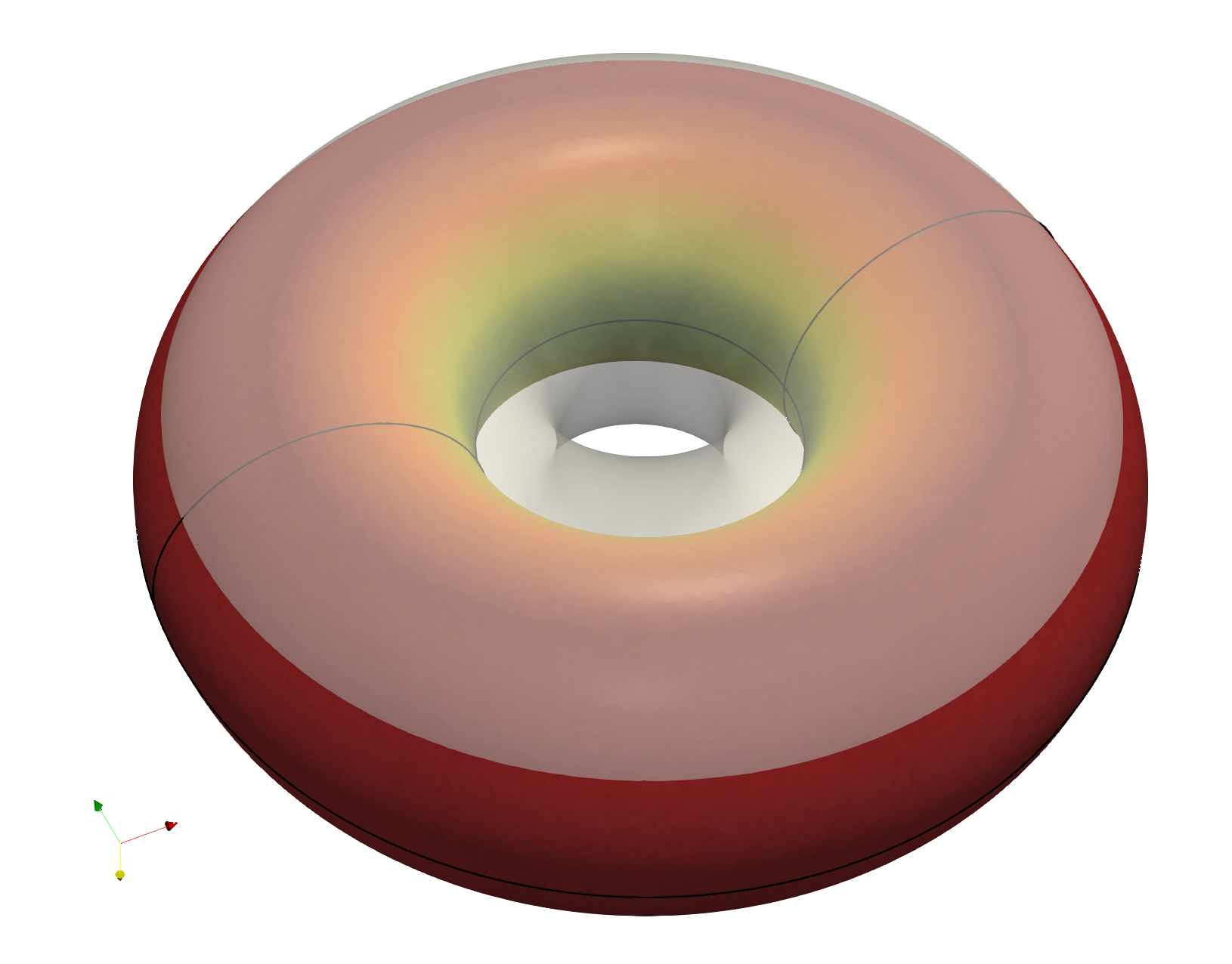}}
        \caption{$t = 0$}
    \end{subfigure}%
    \begin{subfigure}{0.33\textwidth}
        \centering
        \includegraphics[width=\textwidth]{{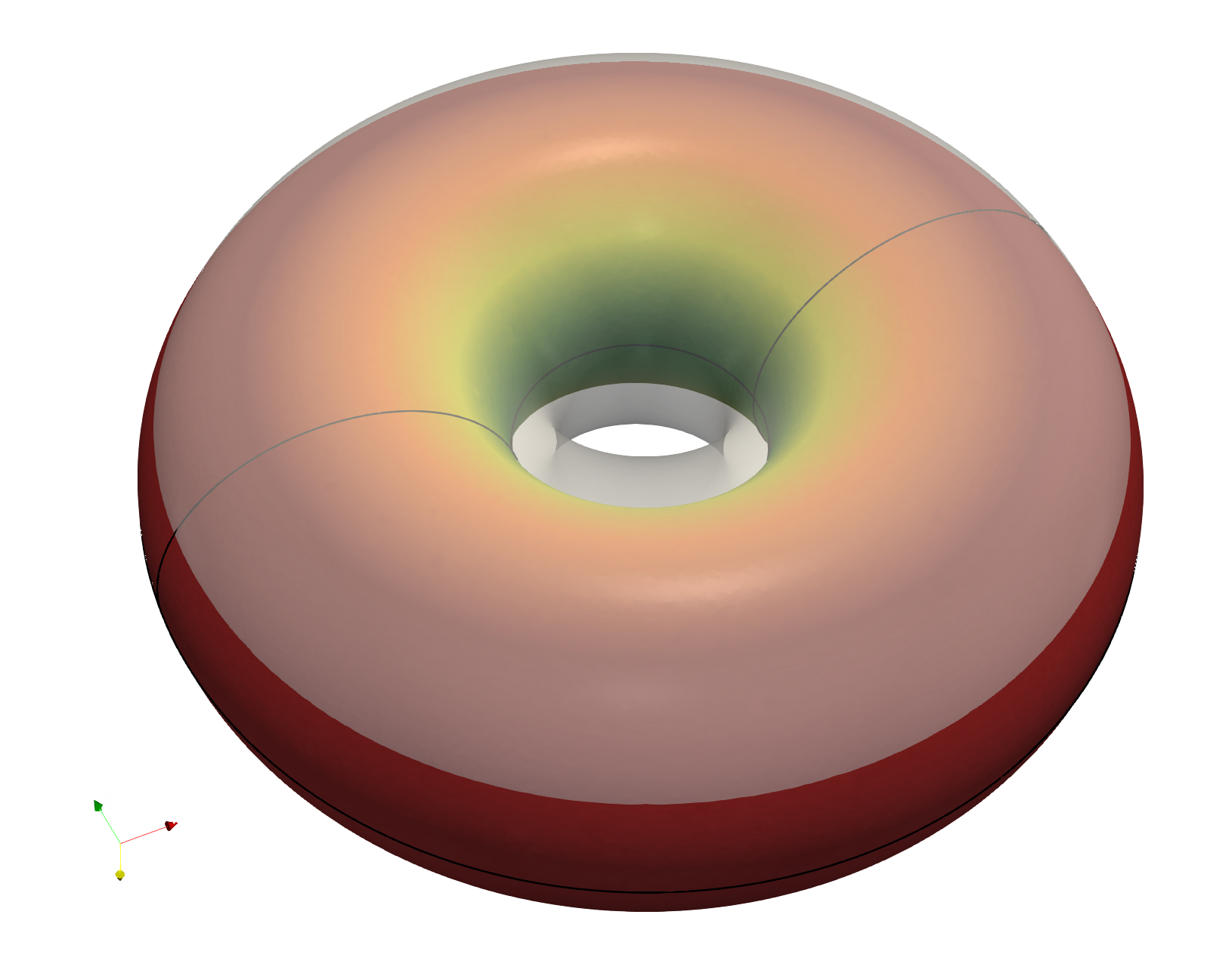}}
        \caption{$t = 0.015$}
    \end{subfigure}%
    \begin{subfigure}{0.33\textwidth}
        \centering
        \includegraphics[width=\textwidth]{{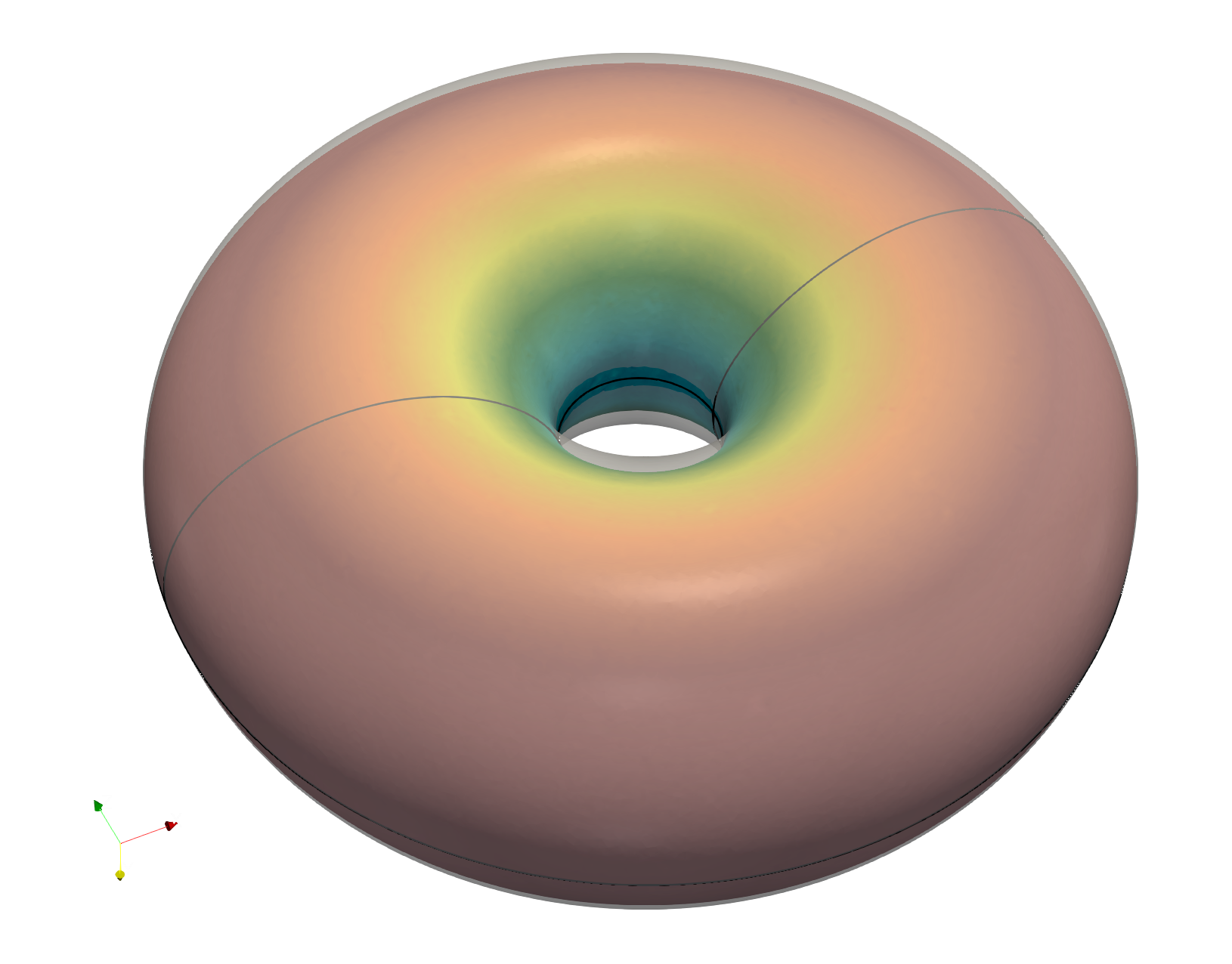}}
        \caption{$t = 0.03$}
    \end{subfigure}%
    \\
    \begin{center}
        \input{figures/willmoretorus3d_eqbshape_colorbar.tex}%
    \end{center}
    \caption{
        3D Willmore flow of a torus toward the Clifford torus.
        The evolving surface is colored with the value of the curvature while
        the Clifford torus is represented in semi-transparent white.
    }
    \label{fig:WillmoreTorus3DEvol}
\end{figure}

\begin{figure}[!htb]
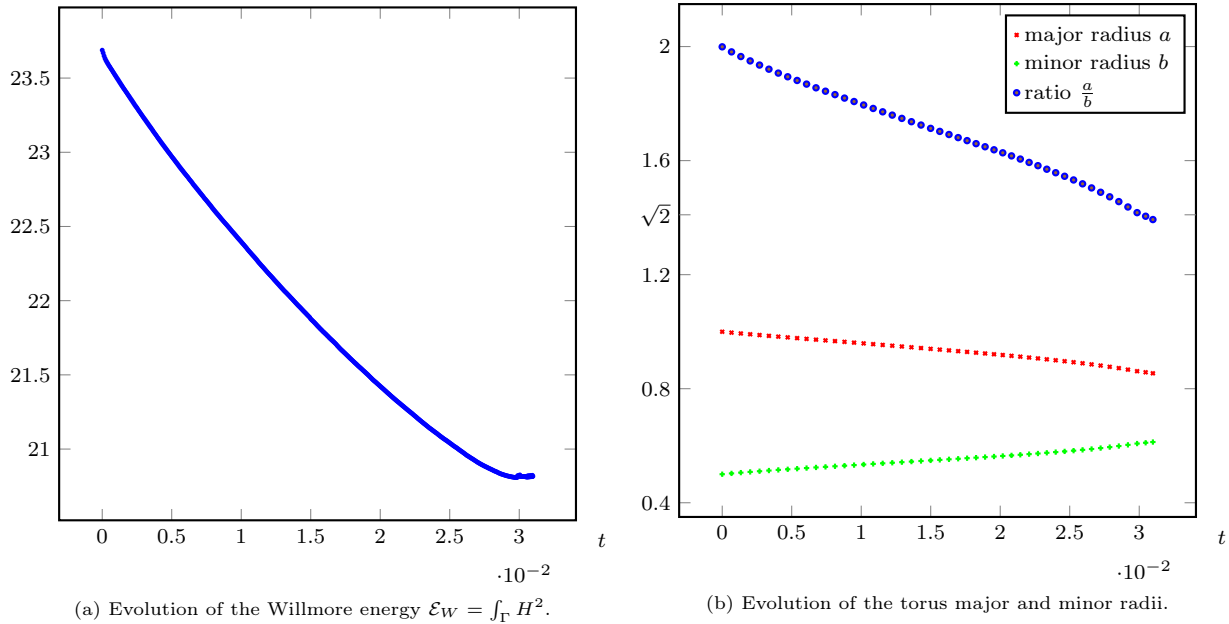

    \centering
    \begin{subfigure}{0.5\textwidth}
        \centering
        \input{figures/willmoretorus3d_lambda2_energy.tex}
        \caption{
            Evolution of the Willmore energy $\mathcal{E}_{W} = \int_{\Gamma} H^2$.
        }
    \end{subfigure}%
    \begin{subfigure}{0.5\textwidth}
        \centering
        \input{figures/willmoretorus3d_lambda2_radii.tex}
        \caption{
            Evolution of the torus major and minor radii.
        }
    \end{subfigure}%
    \caption{
        Willmore energy and geometrical quantities of the torus. We observe
        that the flow decreases the energy to about $20.8$, which is within $5\%$
        from the theoretical Clifford torus energy $2\pi^2 \approx 19.74$. We also
        observe that at equilibrium, the ratio between the radii is indeed about
        $\sqrt{2}$. \review{Beyond $t=3$, the surface leaves the toroidal class, and the computed $a$ and $b$ are not anymore relevant (see above).}
    }
    \label{fig:WillmoreTorus3dMeasures}
\end{figure}

\Cref{fig:WillmoreTorus3DEvol,fig:WillmoreTorus3dMeasures} show the results of
our three-dimensional Willmore flow algorithm \labelcref{alg:3dWillmoreFlow}
starting from a torus with major $a=2$ and minor $b=0.5$ radii. The simulation
was run in a cuboid with lengths $6 \times 4 \times 6$ and an unstructured mesh
of typical size $h \approx 0.02$.  We observe that the surface flows to the
expected Clifford torus, and that the estimated resulting energy and radii ratio
are in good agreement with the expected values for such a torus. \review{However
we were unable to stabilise the optimal shape once it reaches the Clifford
torus. We observe in \Cref{fig:WillmoreTorus3DEvol} the energy fluctuating and
even increasing slowly beyond $t=3$, while in \Cref{fig:WillmoreTorus3dMeasures}
the ratio $a/b$ does not seem to converge. The reason behind this behaviour is
that we are computing $a$ and $b$  from the surface and enclosed volume assuming
that we are in the class of toroidal surfaces. However it turns out that our
algorithm, when reaching the Clifford torus, tries to minimise the energy by
slipping out of this class (we remind that we are making a full 3D computation
not imposing any symmetry). There the central and peripheral circles have no
longer exactly the same center. It seems that while reaching the Clifford torus,
the energy is flattening and it is more likely that discretisation errors drive
the surface out of the tori class. Once this symmetry breaking has occurred, the
hole is then quickly closed and we fall in the $0$ genus surfaces class of lower
energies. This raises the question on how to constraint somehow the surface to
stay in the tori class (without doing axisymmetric computations). For now we did
not find any simple algorithm to tackle that problem, which while interesting,
is not directly related to our aimed application to red blood cells.}

\subsection{Vesicles equilibrium shapes}

We now consider the simulation of vesicles, which provide the most popular and
simple model for capsules, closed thin shell or cells, like human red blood
cells. From a physical point of view, the model consists in an elastic thin
shell of fixed area enclosing some incompressible fluid, so that the inner
volume is also fixed. As such, the vesicle only deforms through bending of the
shell, and is thus controlled by a Willmore-like -- or Canham-Helfrich energy
\cite{helfrich1973elastic,canham1970minimum} (see e.g.
    \cite{seifert1997configurations} for a presentation of the model in a
biological context), which reduces to the standard Willmore energy when the
membrane has no spontaneous curvature \reviewII{-- the unconstrained (open) membrane is flat}, 
which we assume in this work.

In the following, we shall use our diffusion-redistanciation algorithms with
the constraints of constant volume and area of the evolving surface to compute the
2D and 3D equilibrium shapes of vesicles with various reduced volumes. 
\review{We show in the following the results of the ``Euler'' variants of
\cref{alg:2dWillmoreFlowWithConservations} and
\cref{alg:3dWillmoreFlowWithConservations}. We have also performed the same
simulations using the ``Crank-Nicolson'' variants of the algorithms, which
yielded very similar shapes and results, but featured some small -- stable -- oscillations at
the equilibrium, as could be expected from a non-diffusive discretisation scheme.}
The resulting equilibrium shapes are compared to the ones from the literature when
available, and assessed from energetical and stability points of view to
illustrate the good properties of our numerical method.

\subsubsection{Computation of the equilibrium shapes of 2D vesicles}

We first consider the two-dimensional case, and compute the Willmore flow under
the constraint of constant inner volume and surface area of an initial ellipse
with semi-minor and -major axes adjusted according to a chosen reduced volume
\begin{equation}
    \nu = \frac{4\pi V}{A^2 } = \frac{V}{\pi  \left[ \frac{A}{2\pi} \right]^2}
\end{equation}
The reduced volume represents the ratio between the volume of the ellipse and
the volume of a circle of same area.  Note that the conservation of both the
volume and area of the vesicle naturally entails that the reduced volume is also
preserved.

\begin{figure}[!htb]
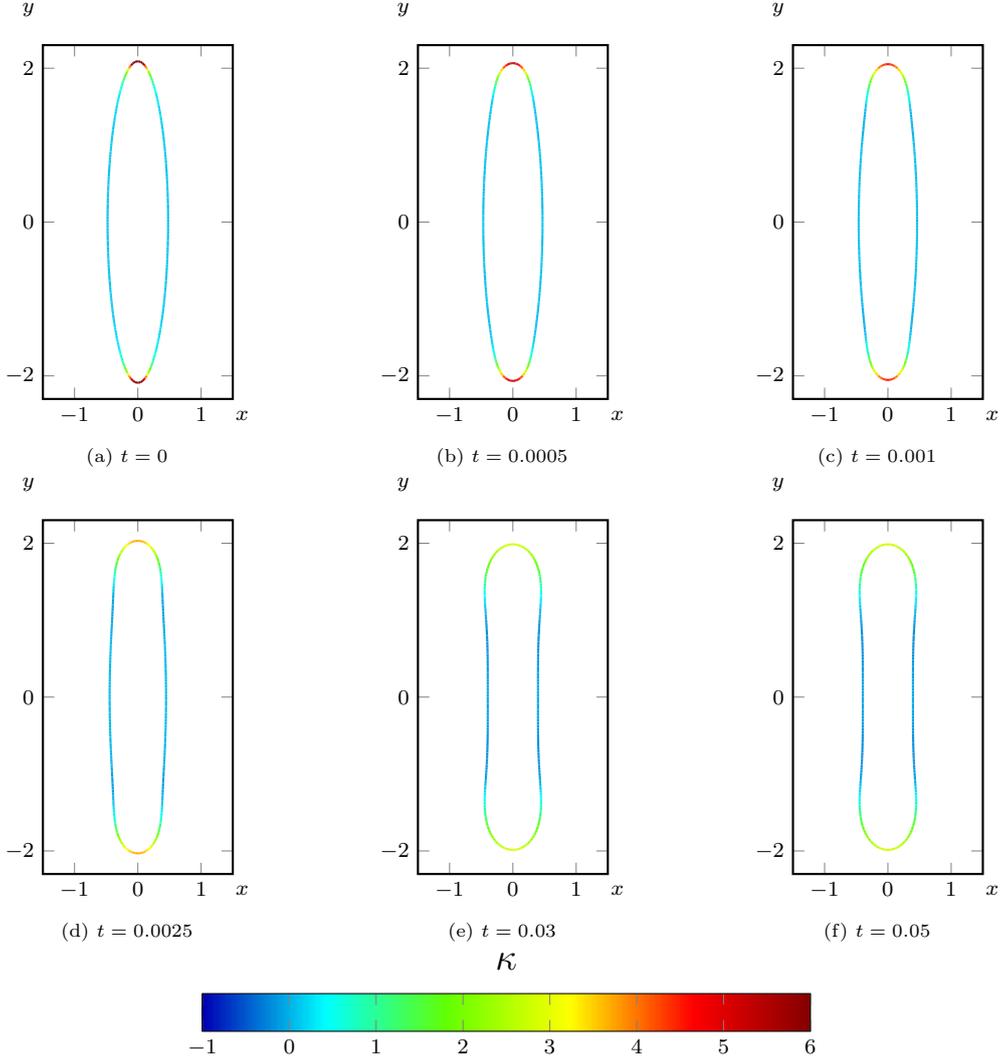

    \centering
    \begin{subfigure}{0.3\textwidth}
        \centering
        \input{figures/vesicle2d_alpha0.5_t0.0.tex}%
        \caption{$t=0$}
    \end{subfigure}%
    \begin{subfigure}{0.3\textwidth}
        \centering
        \input{figures/vesicle2d_alpha0.5_t0.0005.tex}%
        \caption{$t=0.0005$}
    \end{subfigure}%
    \begin{subfigure}{0.3\textwidth}
        \centering
        \input{figures/vesicle2d_alpha0.5_t0.001.tex}%
        \caption{$t=0.001$}
    \end{subfigure}%
    \\
    \begin{subfigure}{0.3\textwidth}
        \centering
        \input{figures/vesicle2d_alpha0.5_t0.0025.tex}%
        \caption{$t=0.0025$}
    \end{subfigure}%
    \begin{subfigure}{0.3\textwidth}
        \centering
        \input{figures/vesicle2d_alpha0.5_t0.03.tex}%
        \caption{$t=0.03$}
    \end{subfigure}%
    \begin{subfigure}{0.3\textwidth}
        \centering
        \input{figures/vesicle2d_alpha0.5_t0.05.tex}%
        \caption{$t=0.05$}
    \end{subfigure}%
    \\
    \ref{2dalpha05itercolorbar}
    \caption{
        Convergence of \cref{alg:2dWillmoreFlowWithConservations} to an
        equilibrium shape for the $\nu=0.5$ case, starting from an ellipse
        with semi-major and -minor axes $2.09021$ and $0.478421$.
        \reviewII{The curves are colored with the value of the curvature.}
    }
    \label{s5:fig:EvolutionToEquilibrium2D}
\end{figure}

\begin{figure}[!ht]
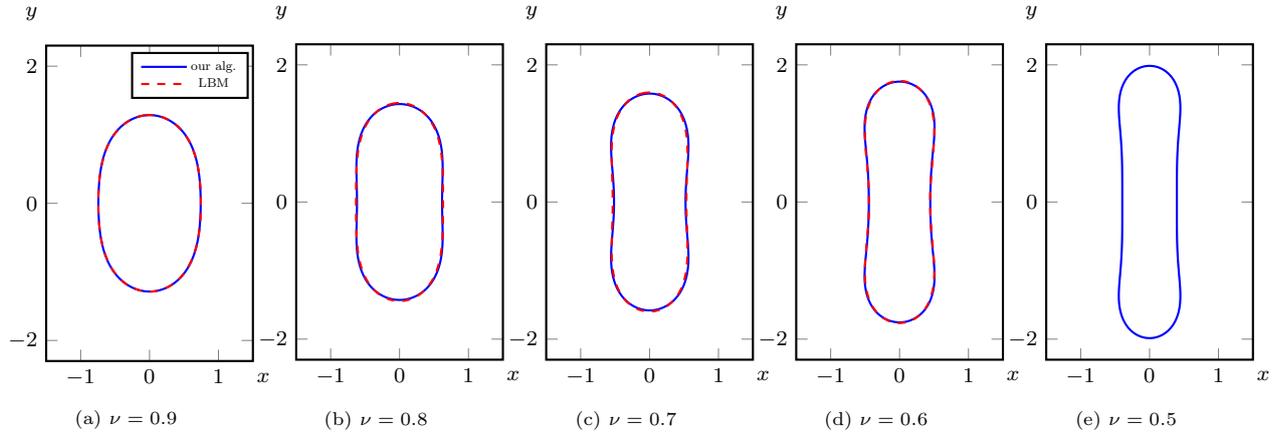

    \centering
    \begin{subfigure}{0.2\textwidth}
        \centering
        \input{figures/vesicle2d_alpha0.9_cmplbm.tex}%
        \caption{$\nu = 0.9$}
    \end{subfigure}%
    \begin{subfigure}{0.2\textwidth}
        \centering
        \input{figures/vesicle2d_alpha0.8_cmplbm.tex}
        \caption{$\nu = 0.8$}
    \end{subfigure}%
    \begin{subfigure}{0.2\textwidth}
        \centering
        \input{figures/vesicle2d_alpha0.7_cmplbm.tex}
        \caption{$\nu = 0.7$}
    \end{subfigure}%
    \begin{subfigure}{0.2\textwidth}
        \centering
        \input{figures/vesicle2d_alpha0.6_cmplbm.tex}
        \caption{$\nu = 0.6$}
    \end{subfigure}%
    \begin{subfigure}{0.2\textwidth}
        \centering
        \input{figures/vesicle2d_alpha0.5_cmplbm.tex}
        \caption{$\nu = 0.5$}
    \end{subfigure}%
    \caption{
        Comparison of the 2D equilibrium shapes obtained by
        \cref{alg:2dWillmoreFlowWithConservations} with the ones obtained by
        \cite{kaoui2011two} using a lattice-Boltzmann method for reduced volumes
        $\nu=0.9$, $0.8$, $0.7$ and $0.6$. 
        We also give our result for a reduced volume $\nu=0.5$.
    }
    \label{s5:fig:ShapesRemplissages2D}
\end{figure}

\begin{figure}[!ht]
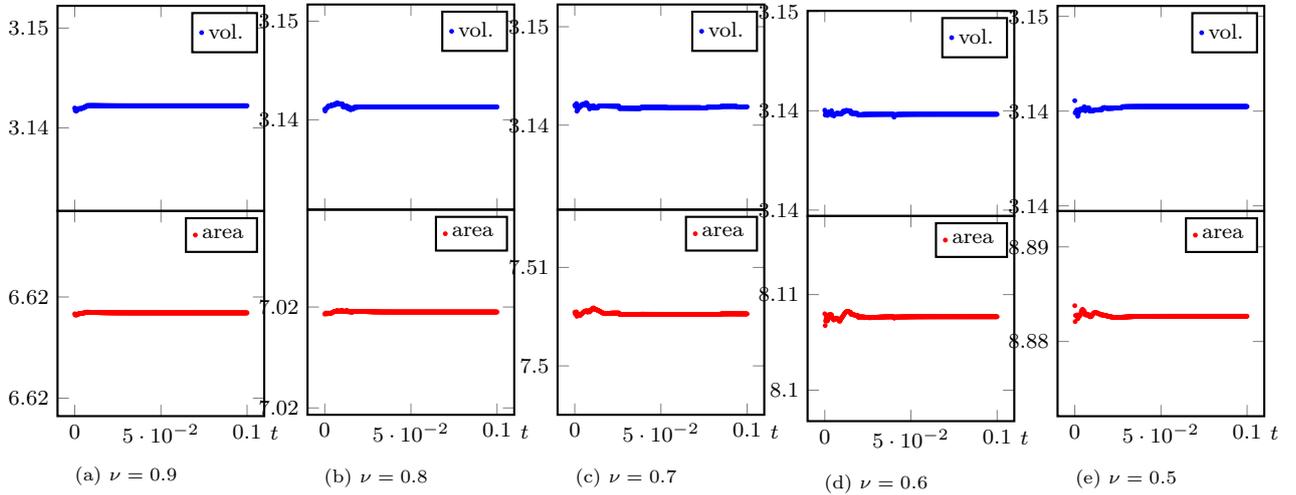

    \centering
    \begin{subfigure}{0.2\textwidth}
        \centering
        \input{figures/vesicle2d_alpha0.9_volarea.tex}%
        \caption{$\nu = 0.9$}
    \end{subfigure}%
    \begin{subfigure}{0.2\textwidth}
        \centering
        \input{figures/vesicle2d_alpha0.8_volarea.tex}
        \caption{$\nu = 0.8$}
    \end{subfigure}%
    \begin{subfigure}{0.2\textwidth}
        \centering
        \input{figures/vesicle2d_alpha0.7_volarea.tex}
        \caption{$\nu = 0.7$}
    \end{subfigure}%
    \begin{subfigure}{0.2\textwidth}
        \centering
        \input{figures/vesicle2d_alpha0.6_volarea.tex}
        \caption{$\nu = 0.6$}
    \end{subfigure}%
    \begin{subfigure}{0.2\textwidth}
        \centering
        \input{figures/vesicle2d_alpha0.5_volarea.tex}
        \caption{$\nu = 0.5$}
    \end{subfigure}%
    \caption{
        Evolution of the 2D vesicle volume and ``area'' (perimeter in this case)
        with the volume- and area-preserving \cref{alg:2dWillmoreFlowWithConservations}. 
        The algorithm displays very good conservation properties : the relative volume
        and area changes are less than respectively $1\cdot 10^{-4}$ and
        $5 \cdot 10^{-4}$ for the most
        difficult $\nu=0.5$ case, and the final equilibrium shapes volumes and
        areas relative changes are $\lesssim 1\cdot 10^{-5}$ for all the cases.
    }
    \label{s5:fig:vesicle2dVolArea}
\end{figure}

\begin{figure}[!ht]
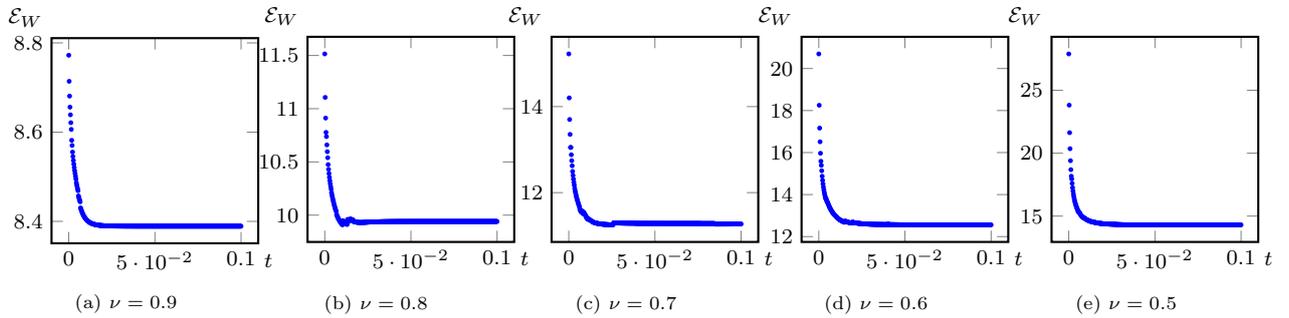

    \centering
    \begin{subfigure}{0.2\textwidth}
        \centering
        \input{figures/vesicle2d_alpha0.9_energy.tex}%
        \caption{$\nu = 0.9$}
    \end{subfigure}%
    \begin{subfigure}{0.2\textwidth}
        \centering
        \input{figures/vesicle2d_alpha0.8_energy.tex}
        \caption{$\nu = 0.8$}
    \end{subfigure}%
    \begin{subfigure}{0.2\textwidth}
        \centering
        \input{figures/vesicle2d_alpha0.7_energy.tex}
        \caption{$\nu = 0.7$}
    \end{subfigure}%
    \begin{subfigure}{0.2\textwidth}
        \centering
        \input{figures/vesicle2d_alpha0.6_energy.tex}
        \caption{$\nu = 0.6$}
    \end{subfigure}%
    \begin{subfigure}{0.2\textwidth}
        \centering
        \input{figures/vesicle2d_alpha0.5_energy.tex}
        \caption{$\nu = 0.5$}
    \end{subfigure}%
    \caption{
        Evolution of the Willmore energy $\mathcal{E}_{W} = \int_{\Gamma} H^2$ 
        \labelcref{eq:willmore_energy} of the 2D vesicle 
        with the volume- and area-preserving 
        \cref{alg:2dWillmoreFlowWithConservations}. 
    }
    \label{s5:fig:vesicle2dEnergy}
\end{figure}

The simulations were run in a square domain $[-4,4]^2$ with a structured
triangle mesh with mesh size $h = 0.04$ and adaptive time-step strategy. We
illustrate the convergence of our method to an equilibrium shape in
\cref{s5:fig:EvolutionToEquilibrium2D} and compare our results with the ones
obtained using a lattice-Boltzmann method in \cite{kaoui2011two} in
\cref{s5:fig:ShapesRemplissages2D}. We observe that the shapes computed with our
algorithm are in excellent agreement with the ones obtained from direct
lattice-Boltzmann simulations.  We also illustrate the good volume and area
conservation properties of our algorithm in \cref{s5:fig:vesicle2dVolArea}. The
evolution of the Willmore energy along the flow is shown in
\cref{s5:fig:vesicle2dEnergy}.

\subsubsection{Computation of vesicles equilibrium shapes in 3D}

The landscape of vesicle equilibrium shapes is much more rich in three
dimensions, as the Gauss curvature comes into play and allows for multiple local
minima for shapes with the same topological class and invariants.  Axisymmetric
equilibrium shapes have been mapped out in a phase diagram in
\cite{seifert1991shape}, and some non-axisymmetric shapes have been computed in
\cite{julicher1993conformal,du2006simulating,du2008adaptive,bretin2015phase},
but the general phase diagram for Willmore energy minimisers is still an open
question.

In our case, we consider the evolution of ellipsoids under the constant volume
and area Willmore flow. As such, we expect to recover the prolate and oblate
axisymmetric shapes of zero genus and spontaneous curvature referenced in
\cite{seifert1991shape} depending on the initial -- conserved -- reduced volume,

\begin{equation} 
    \nu 
    = \frac{6\sqrt{\pi}V}{A^\frac{3}{2}}
    = \frac{V}{\frac43\pi \left[\frac A{4\pi}\right]^\frac32}.
\end{equation}

As shown in \cite{seifert1991shape}, the prolate and oblate shapes are
both local minimisers for $\nu \approx 0.51$. While the prolate shapes exist for all
possible $\nu$, the oblate ones self-intersect below $\nu \approx 0.51$, and the
stomatocyte shapes become the only non-prolate feasible minimisers.
In this work, we restrict ourselves to $\nu \geq 0.6$, and study both the oblate
and prolate cases, as our algorithm appears able to capture the local minimising
shapes in a stable way.

Fewer quantitative works are available for comparison in this three dimension
case. We can mention \cite{biben2005phase,du2006simulating,
feng2006finite,du2008adaptive} where a mix of qualitative and quantitative
results are presented. In \cite{biben2005phase,du2006simulating}, a phase field
method  was developed and qualitative results were given. In
\cite{du2008adaptive} an adaptive version of the phase field model was
proposed with some quantitative results in terms of energy, which are hard to
compare with theoretical known values of Seifert \cite{seifert1991shape}. The
approach of Feng and Klug \cite{feng2006finite} is based on surface finite
element and is more quantitative in its results as far as the energy is concerned. 

\paragraph{Oblate case}

The oblate case simulations were run in a cuboid with lengths $3.6 \times
6.4 \times 6.4$ and mesh size $h \approx 0.04$.
The initial level-set functions were taken as signed distance functions to oblate
ellipsoids with semi-minor axis $a$ and semi-major axes $b = c$ adjusted to
get a volume $\frac{4\pi}{3}$ for all the simulations and a corresponding
reduced volume $\nu = 0.6$, $0.65$, $0.7$, $0.8$ or $0.9$.

\begin{figure}[!htb]
    \centering
    \begin{subfigure}{0.2\textwidth}
        \centering
        \includegraphics[width=\textwidth]{{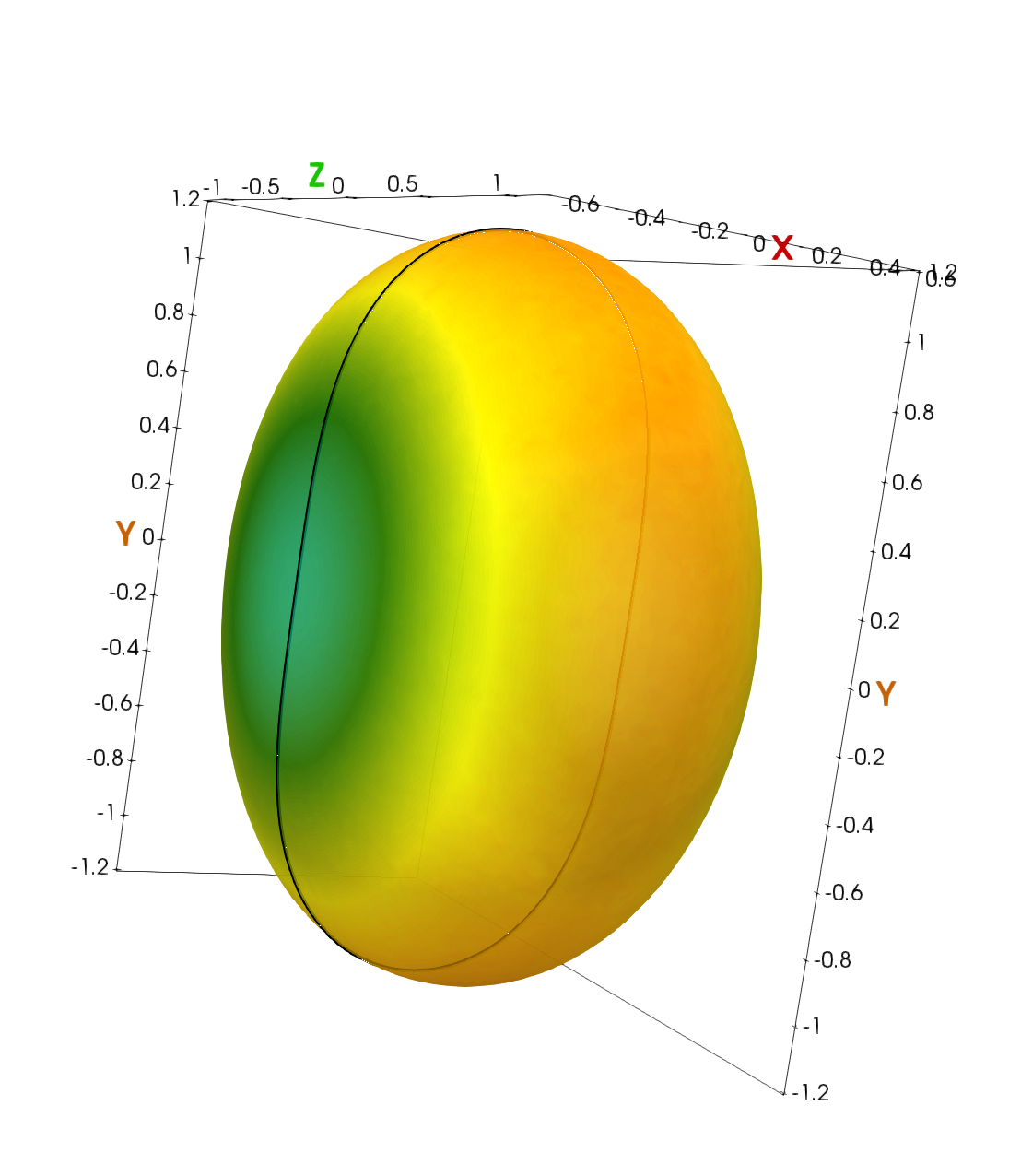}}
        \caption{$\nu = 0.9$}
    \end{subfigure}%
    \begin{subfigure}{0.2\textwidth}
        \centering
        \includegraphics[width=\textwidth]{{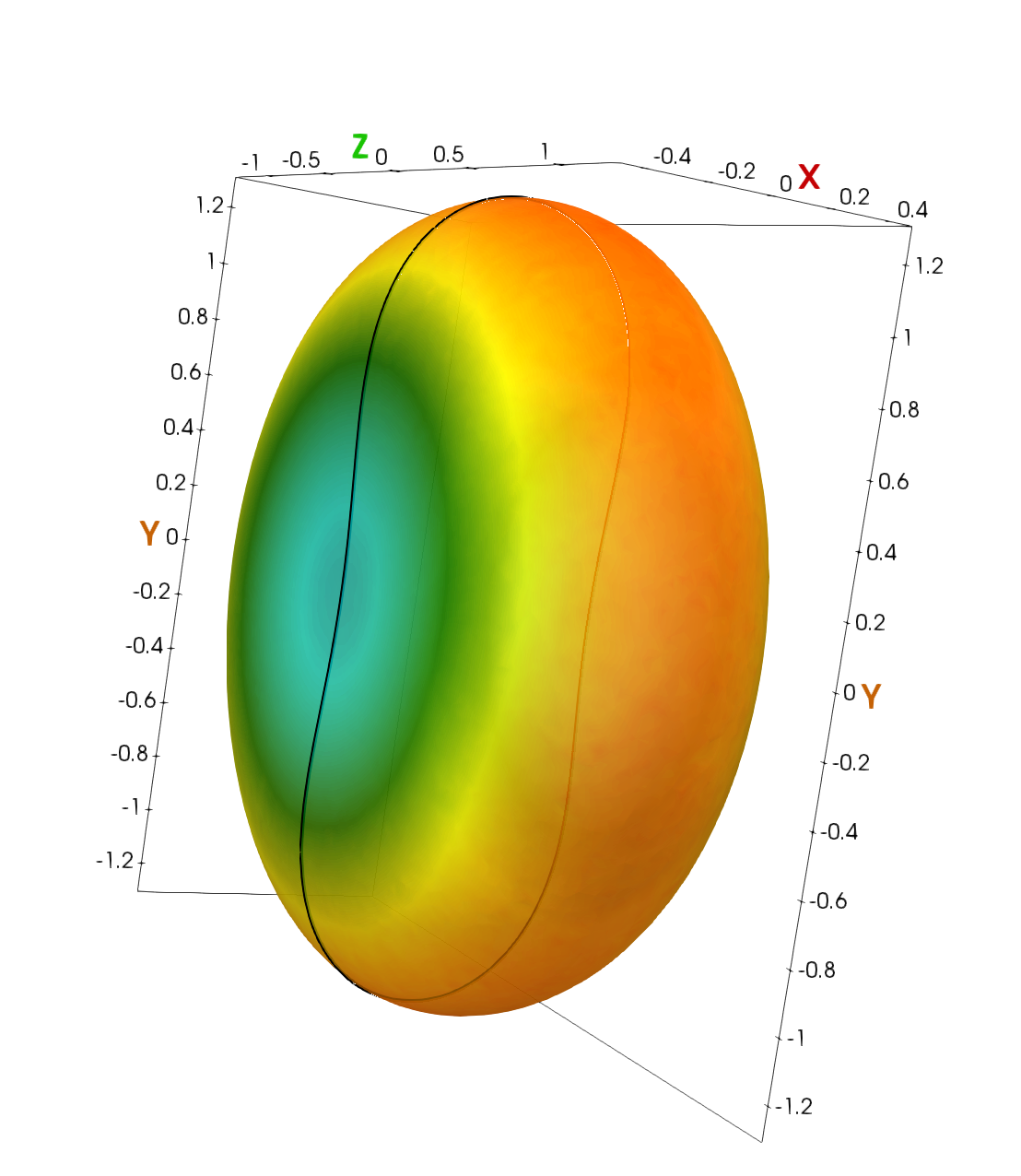}}
        \caption{$\nu = 0.8$}
    \end{subfigure}%
    \begin{subfigure}{0.2\textwidth}
        \centering
        \includegraphics[width=\textwidth]{{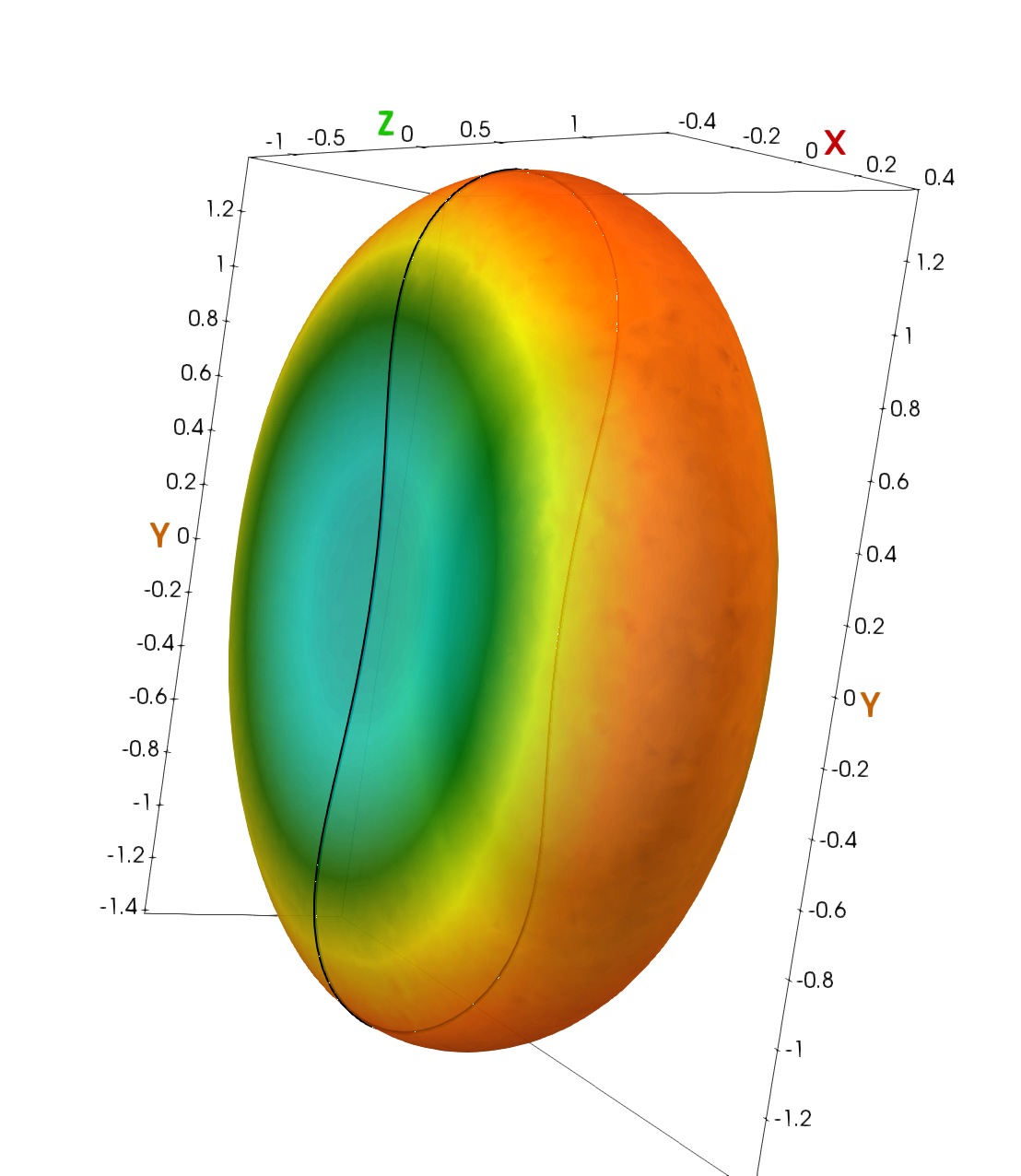}}
        \caption{$\nu = 0.7$}
    \end{subfigure}%
    \begin{subfigure}{0.2\textwidth}
        \centering
        \includegraphics[width=\textwidth]{{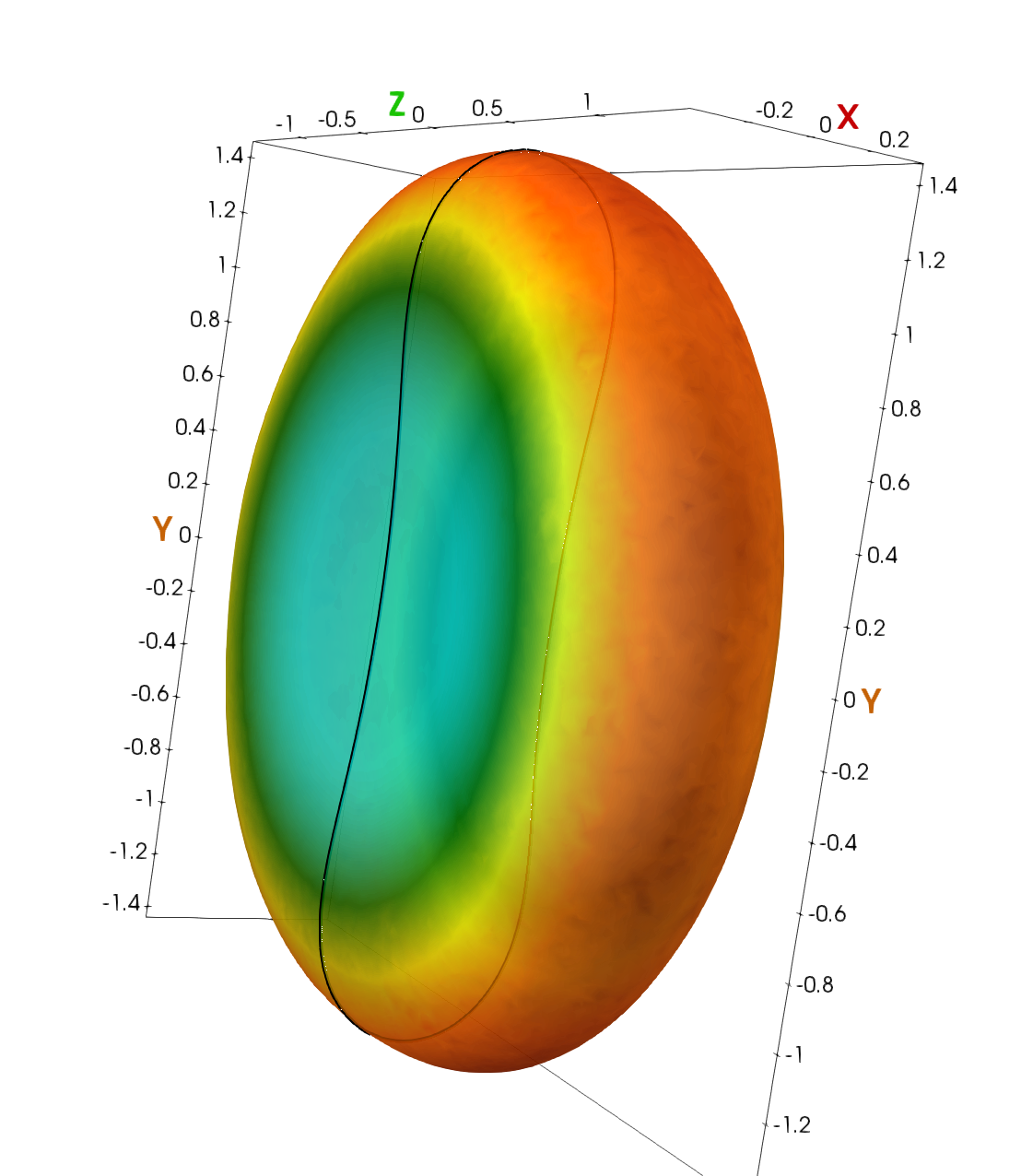}}
        \caption{$\nu = 0.65$}
    \end{subfigure}%
    \begin{subfigure}{0.2\textwidth}
        \centering
        \includegraphics[width=\textwidth]{{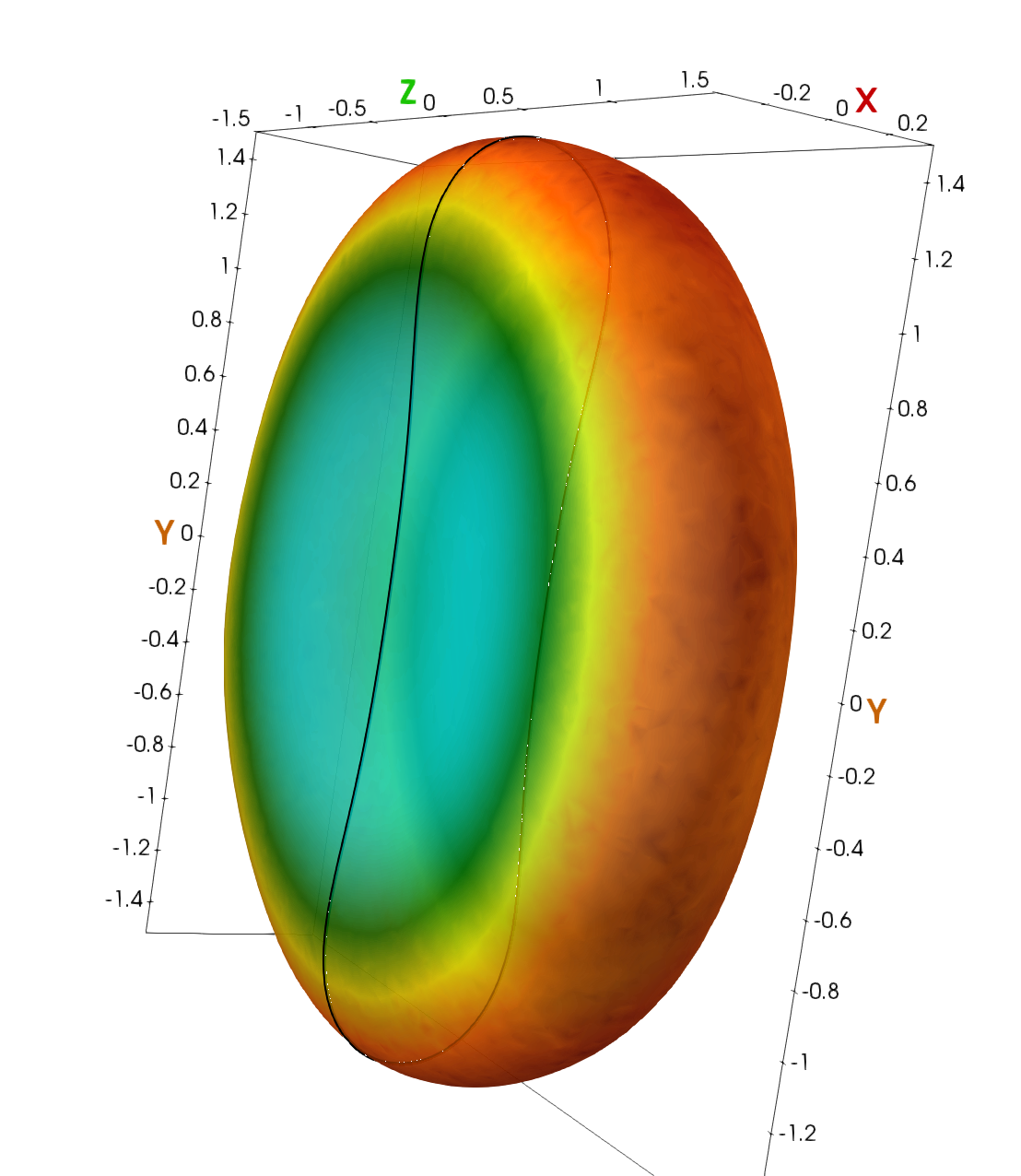}}
        \caption{$\nu = 0.6$}
    \end{subfigure}%
    \\
    \begin{center}
        \input{figures/vesicle3d_eqbshape_colorbar.tex}%
    \end{center}
    \caption{
        3D oblate vesicles equilibrium shapes obtained with
        \cref{alg:3dWillmoreFlowWithConservations} for the reduced volumes
        $\nu=0.9$, $0.8$, $0.7$, $0.65$ and $0.6$.
    }
    \label{sec5:fig:3dOblEqbShapes}
\end{figure}

\Cref{sec5:fig:3dOblEqbShapes} shows the equilibrium shapes obtain with our
constrained diffusion-redistanciation algorithm
\labelcref{alg:3dWillmoreFlowWithConservations} for the
different reduced volumes. We also show two-dimensional cuts of these
equilibrium shapes in \cref{s5:fig:3dOblEqbSlices} to illustrate the good
symmetry conservation properties of our algorithm.
\Cref{s5:fig:vesicle3dOblVolArea,s5:fig:vesicle3dOblEnergy} show the evolution
of the Willmore energy, and the area and enclosed volume of the surface as it
evolves.

\begin{figure}[!ht]
    \centering
    \begin{subfigure}{0.2\textwidth}
        \centering
        \input{figures/vesicle3d_alpha0.9obl_slices.tex}%
        \caption{$\nu = 0.9$}
    \end{subfigure}%
    \begin{subfigure}{0.2\textwidth}
        \centering
        \input{figures/vesicle3d_alpha0.8obl_slices.tex}%
        \caption{$\nu = 0.8$}
    \end{subfigure}%
    \begin{subfigure}{0.2\textwidth}
        \centering
        \input{figures/vesicle3d_alpha0.7obl_slices.tex}%
        \caption{$\nu = 0.7$}
    \end{subfigure}%
    \begin{subfigure}{0.2\textwidth}
        \centering
        \input{figures/vesicle3d_alpha0.65obl_slices.tex}%
        \caption{$\nu = 0.65$}
    \end{subfigure}%
    \begin{subfigure}{0.2\textwidth}
        \centering
        \input{figures/vesicle3d_alpha0.6obl_slices.tex}%
        \caption{$\nu = 0.6$}
    \end{subfigure}%
    \\
    \ref{3doblshapessliceslegend}
    \caption{
        Cut slices in the $x-y$, $x-z$ and $y-z$ planes of the 3D oblate vesicles
        equilibrium shapes displayed in \cref{sec5:fig:3dOblEqbShapes}. These cuts
        highlight the good symmetry preserving property of our algorithm despite the
        absence of any symmetry enforcing method.
    }
    \label{s5:fig:3dOblEqbSlices}
\end{figure}

\begin{figure}[!ht]
    \centering
    \begin{subfigure}{0.2\textwidth}
        \centering
        \input{figures/vesicle3d_alpha0.9obl_volarea.tex}%
        \caption{$\nu = 0.9$}
    \end{subfigure}%
    \begin{subfigure}{0.2\textwidth}
        \centering
        \input{figures/vesicle3d_alpha0.8obl_volarea.tex}
        \caption{$\nu = 0.8$}
    \end{subfigure}%
    \begin{subfigure}{0.2\textwidth}
        \centering
        \input{figures/vesicle3d_alpha0.7obl_volarea.tex}
        \caption{$\nu = 0.7$}
    \end{subfigure}%
    \begin{subfigure}{0.2\textwidth}
        \centering
        \input{figures/vesicle3d_alpha0.65obl_volarea.tex}
        \caption{$\nu = 0.65$}
    \end{subfigure}%
    \begin{subfigure}{0.2\textwidth}
        \centering
        \input{figures/vesicle3d_alpha0.6obl_volarea.tex}
        \caption{$\nu = 0.6$}
    \end{subfigure}%
    \caption{
        Evolution of the 3D oblate vesicle volume and area 
        with the volume- and area-preserving \cref{alg:3dWillmoreFlowWithConservations}. 
        The 3D algorithm also displays good conservation properties: the relative volume
        and area changes are respectively $\lesssim 5\cdot 10^{-5}$ and
        $\lesssim 2 \cdot 10^{-4}$ for all the cases.
    }
    \label{s5:fig:vesicle3dOblVolArea}
\end{figure}

\begin{figure}[!ht]
    \centering
    \begin{subfigure}{0.2\textwidth}
        \centering
        \input{figures/vesicle3d_alpha0.9obl_energy.tex}%
        \caption{$\nu = 0.9$}
    \end{subfigure}%
    \begin{subfigure}{0.2\textwidth}
        \centering
        \input{figures/vesicle3d_alpha0.8obl_energy.tex}
        \caption{$\nu = 0.8$}
    \end{subfigure}%
    \begin{subfigure}{0.2\textwidth}
        \centering
        \input{figures/vesicle3d_alpha0.7obl_energy.tex}
        \caption{$\nu = 0.7$}
    \end{subfigure}%
    \begin{subfigure}{0.2\textwidth}
        \centering
        \input{figures/vesicle3d_alpha0.65obl_energy.tex}
        \caption{$\nu = 0.65$}
    \end{subfigure}%
    \begin{subfigure}{0.2\textwidth}
        \centering
        \input{figures/vesicle3d_alpha0.6obl_energy.tex}
        \caption{$\nu = 0.6$}
    \end{subfigure}%
    \caption{
        Evolution of the Willmore energy $\mathcal{E}_{W} = \int_{\Gamma} H^2$ 
        \labelcref{eq:willmore_energy} of the 3D oblate vesicle 
        with the volume- and area-preserving 
        \cref{alg:3dWillmoreFlowWithConservations}. 
    }
    \label{s5:fig:vesicle3dOblEnergy}
\end{figure}

\paragraph{Prolate case}

The prolate case simulations were run in a cuboid with lengths $9 \times
4 \times 4$ and mesh size $h \approx 0.04$. The initial level-set functions were
chosen as in the oblate case, but using prolate ellipsoids.

\begin{figure}[!htb]
    \centering
    \begin{subfigure}{0.2\textwidth}
        \centering
        \includegraphics[width=\textwidth]{{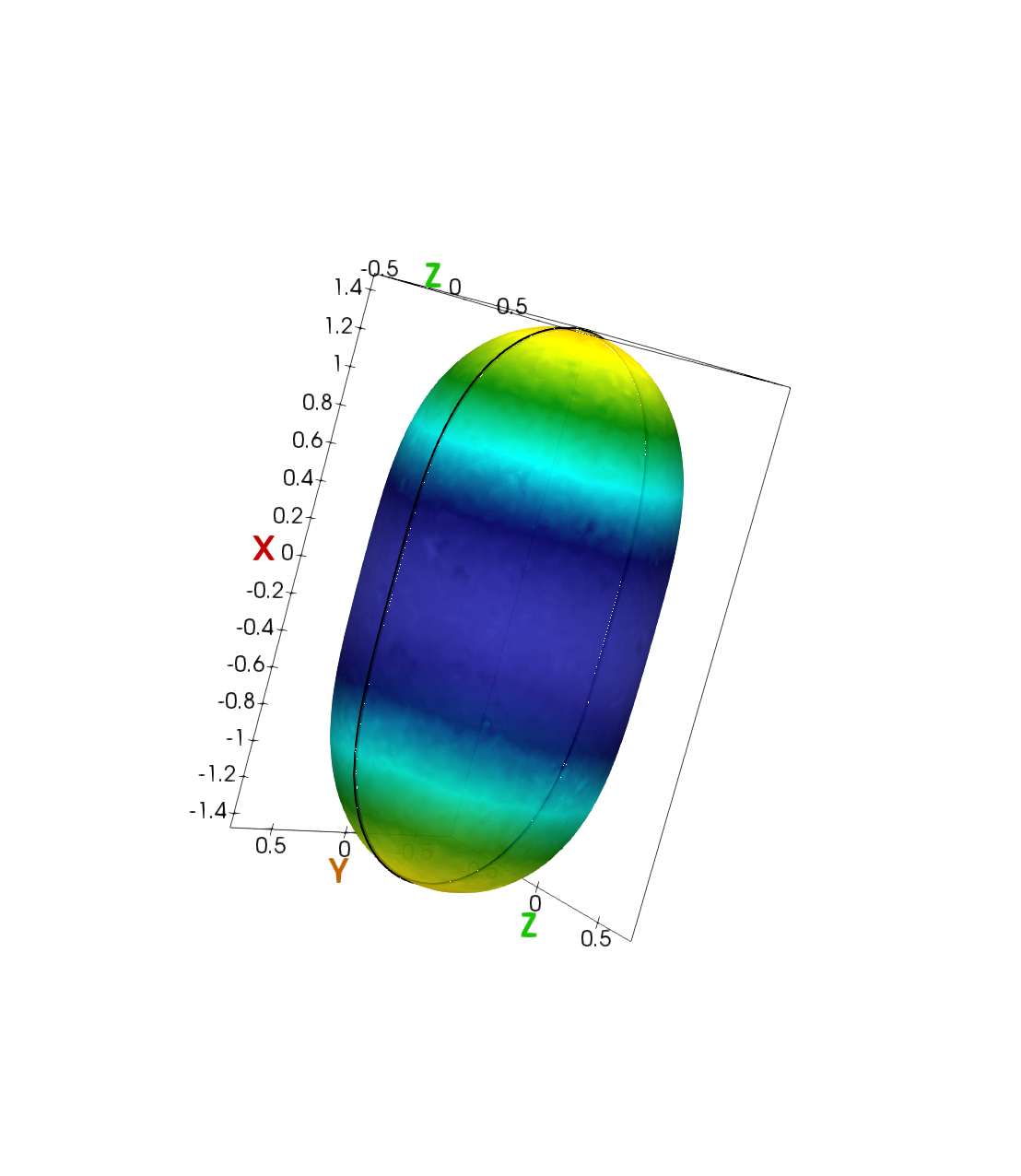}}
        \caption{$\nu = 0.9$}
    \end{subfigure}%
    \begin{subfigure}{0.2\textwidth}
        \centering
        \includegraphics[width=\textwidth]{{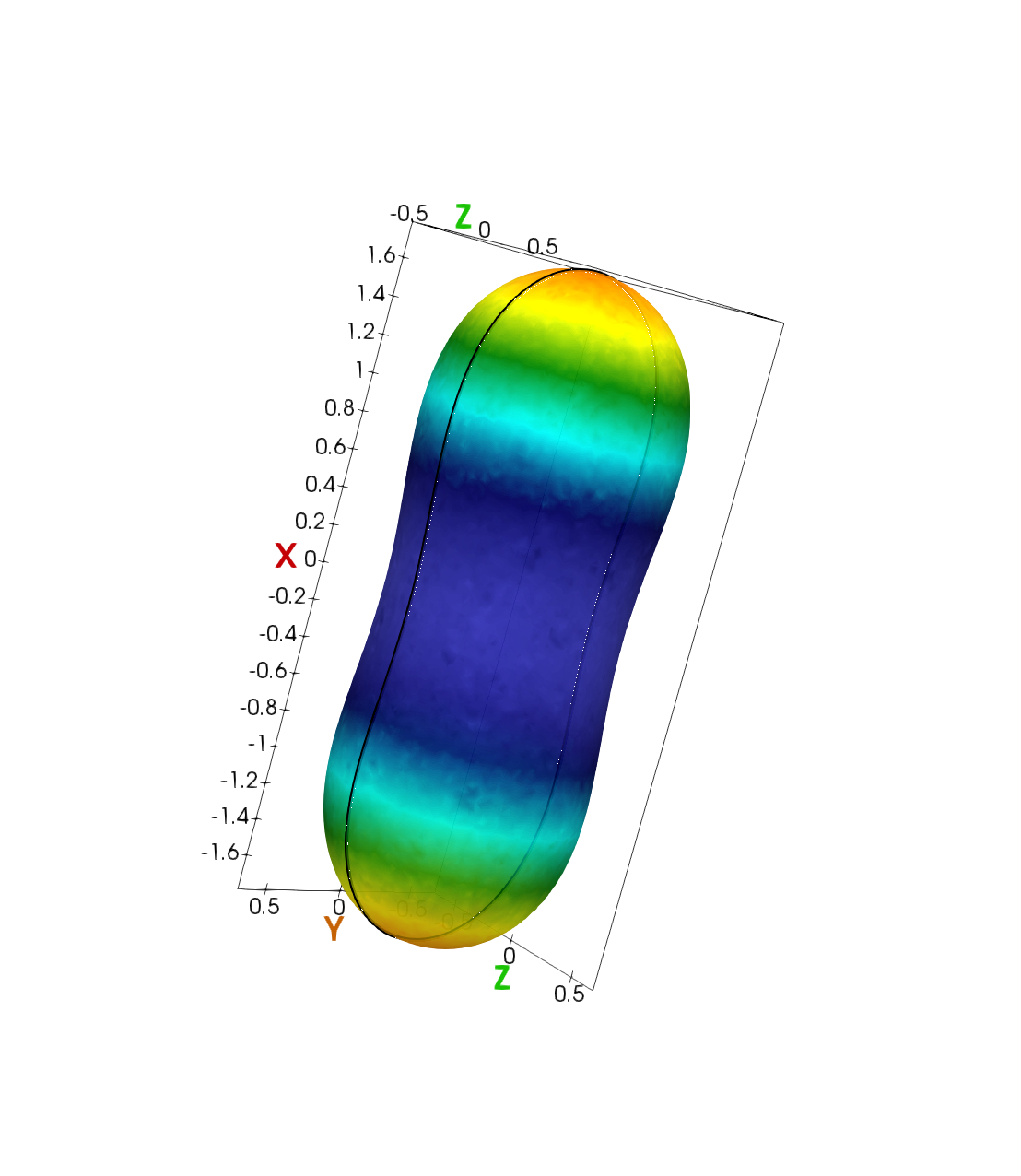}}
        \caption{$\nu = 0.8$}
    \end{subfigure}%
    \begin{subfigure}{0.2\textwidth}
        \centering
        \includegraphics[width=\textwidth]{{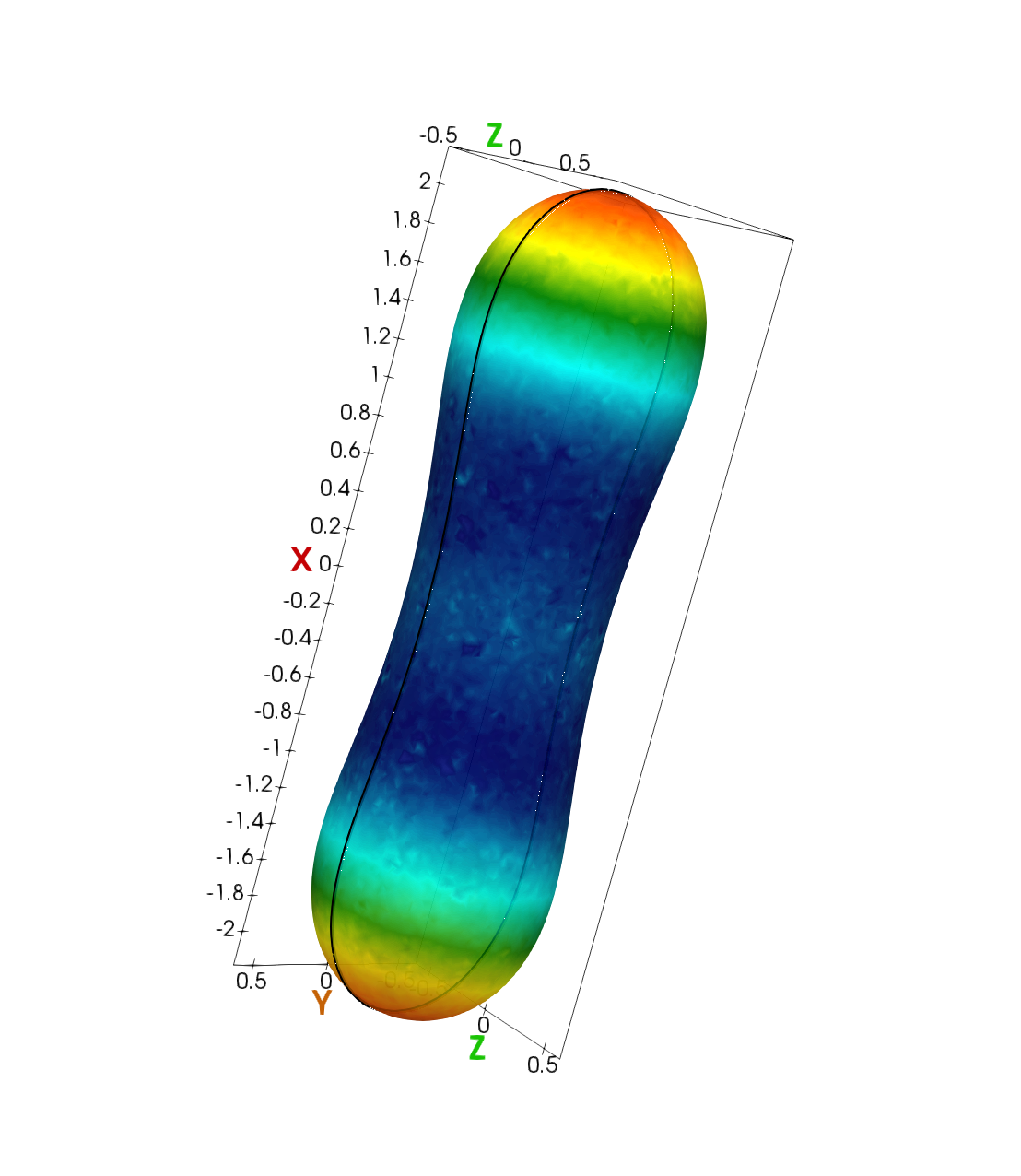}}
        \caption{$\nu = 0.7$}
    \end{subfigure}%
    \begin{subfigure}{0.2\textwidth}
        \centering
        \includegraphics[width=\textwidth]{{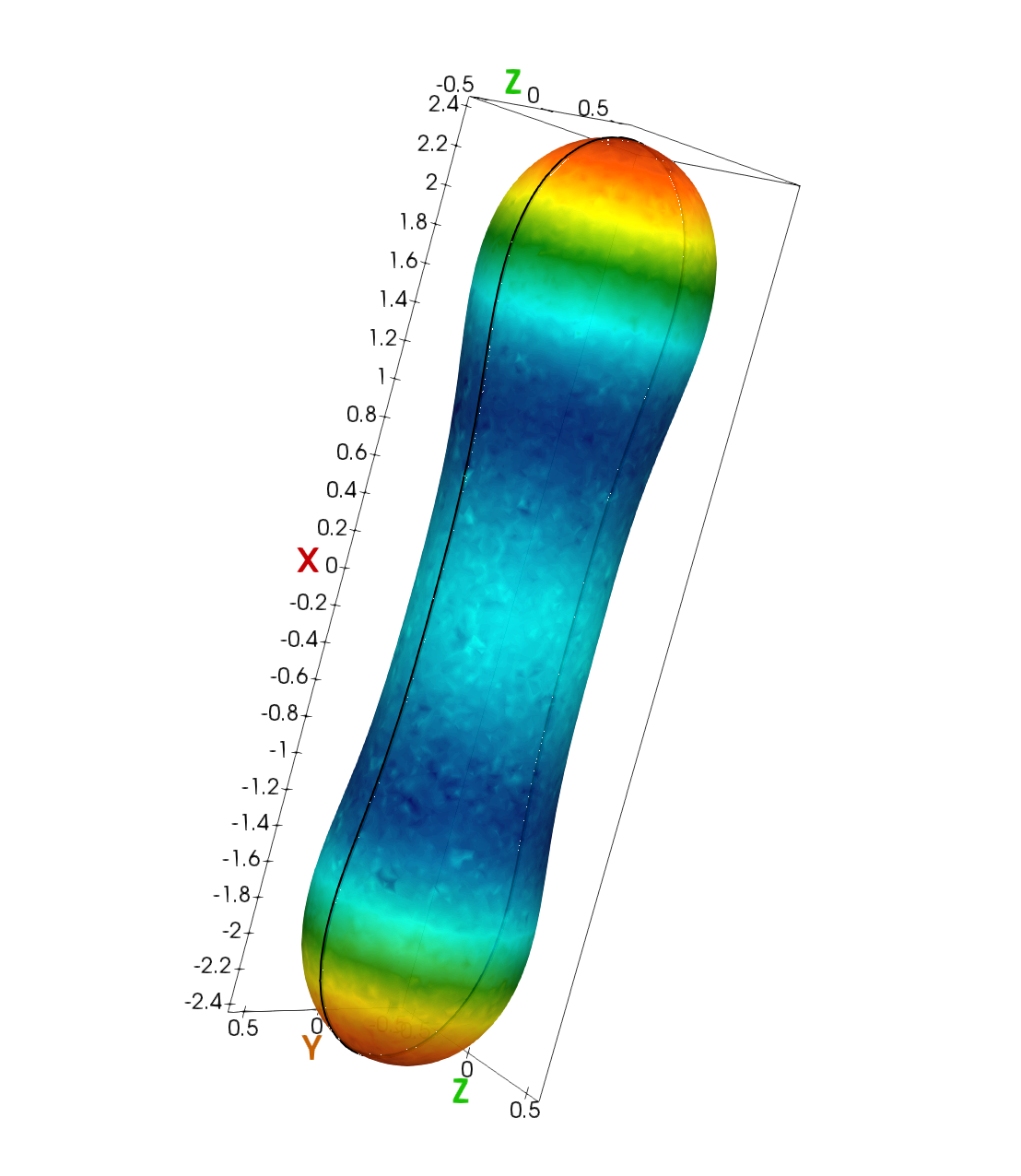}}
        \caption{$\nu = 0.65$}
    \end{subfigure}%
    \begin{subfigure}{0.2\textwidth}
        \centering
        \includegraphics[width=\textwidth]{{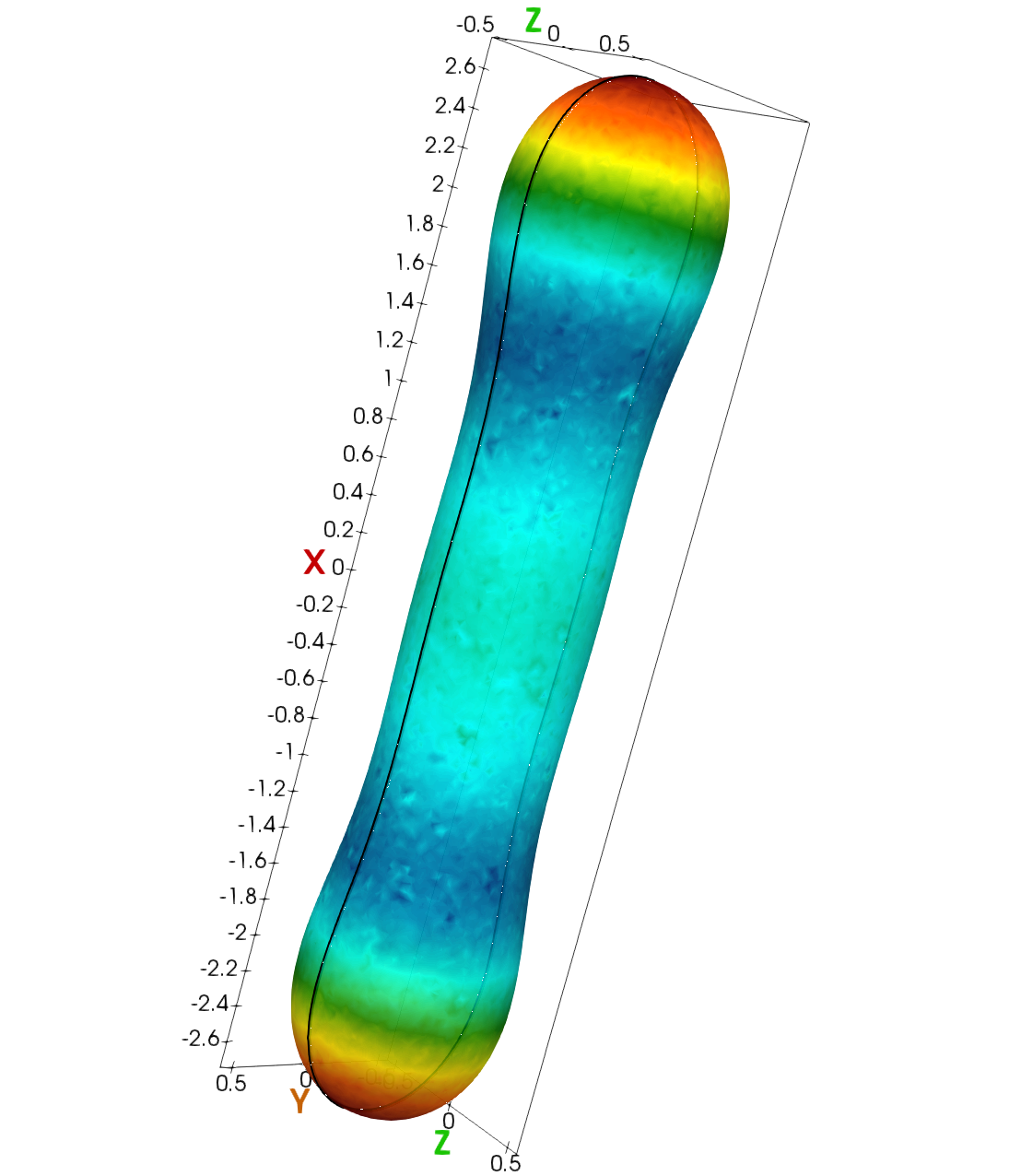}}
        \caption{$\nu = 0.6$}
    \end{subfigure}%
    \\
    \begin{center}
        \input{figures/vesicle3dpro_eqbshape_colorbar.tex}%
    \end{center}
    \caption{
        3D prolate vesicles equilibrium shapes obtained with
        \cref{alg:3dWillmoreFlowWithConservations} for the reduced volumes
        $\nu=0.9$, $0.8$, $0.7$, $0.65$ and $0.6$.
    }
    \label{sec5:fig:3dProEqbShapes}
\end{figure}

The resulting equilibrium shapes and cuts are shown in
\cref{sec5:fig:3dProEqbShapes,s5:fig:3dProEqbSlices} respectively, while
\cref{s5:fig:vesicle3dProVolArea,s5:fig:vesicle3dProEnergy} show the evolution
of the surface energy, area and enclosed volume.

\begin{figure}[!ht]
    \centering
    \begin{subfigure}{0.2\textwidth}
        \centering
        \input{figures/vesicle3d_alpha0.9pro_slices.tex}%
        \caption{$\nu = 0.9$}
    \end{subfigure}%
    \begin{subfigure}{0.2\textwidth}
        \centering
        \input{figures/vesicle3d_alpha0.8pro_slices.tex}%
        \caption{$\nu = 0.8$}
    \end{subfigure}%
    \begin{subfigure}{0.2\textwidth}
        \centering
        \input{figures/vesicle3d_alpha0.7pro_slices.tex}%
        \caption{$\nu = 0.7$}
    \end{subfigure}%
    \begin{subfigure}{0.2\textwidth}
        \centering
        \input{figures/vesicle3d_alpha0.65pro_slices.tex}%
        \caption{$\nu = 0.65$}
    \end{subfigure}%
    \begin{subfigure}{0.2\textwidth}
        \centering
        \input{figures/vesicle3d_alpha0.6pro_slices.tex}%
        \caption{$\nu = 0.6$}
    \end{subfigure}%
    \\
    \ref{3dproshapessliceslegend}
    \caption{
        Cut slices in the $x-y$, $x-z$ and $y-z$ planes of the 3D prolate vesicles
        equilibrium shapes displayed in \cref{sec5:fig:3dOblEqbShapes}. These cuts
        highlight the good symmetry preserving property of our algorithm despite the
        absence of any symmetry enforcing method.
    }
    \label{s5:fig:3dProEqbSlices}
\end{figure}

\begin{figure}[!ht]
    \centering
    \begin{subfigure}{0.2\textwidth}
        \centering
        \input{figures/vesicle3d_alpha0.9pro_volarea.tex}%
        \caption{$\nu = 0.9$}
    \end{subfigure}%
    \begin{subfigure}{0.2\textwidth}
        \centering
        \input{figures/vesicle3d_alpha0.8pro_volarea.tex}
        \caption{$\nu = 0.8$}
    \end{subfigure}%
    \begin{subfigure}{0.2\textwidth}
        \centering
        \input{figures/vesicle3d_alpha0.7pro_volarea.tex}
        \caption{$\nu = 0.7$}
    \end{subfigure}%
    \begin{subfigure}{0.2\textwidth}
        \centering
        \input{figures/vesicle3d_alpha0.65pro_volarea.tex}
        \caption{$\nu = 0.65$}
    \end{subfigure}%
    \begin{subfigure}{0.2\textwidth}
        \centering
        \input{figures/vesicle3d_alpha0.6pro_volarea.tex}
        \caption{$\nu = 0.6$}
    \end{subfigure}%
    \caption{
        Evolution of the 3D prolate vesicle volume and area 
        with the volume- and area-preserving \cref{alg:3dWillmoreFlowWithConservations}. 
        The 3D algorithm also displays good conservation properties: the relative volume
        and area changes are respectively $\lesssim 5\cdot 10^{-6}$ and
        $\lesssim 3 \cdot 10^{-5}$ for all the cases.
    }
    \label{s5:fig:vesicle3dProVolArea}
\end{figure}

\begin{figure}[!ht]
    \centering
    \begin{subfigure}{0.2\textwidth}
        \centering
        \input{figures/vesicle3d_alpha0.9pro_energy.tex}%
        \caption{$\nu = 0.9$}
    \end{subfigure}%
    \begin{subfigure}{0.2\textwidth}
        \centering
        \input{figures/vesicle3d_alpha0.8pro_energy.tex}
        \caption{$\nu = 0.8$}
    \end{subfigure}%
    \begin{subfigure}{0.2\textwidth}
        \centering
        \input{figures/vesicle3d_alpha0.7pro_energy.tex}
        \caption{$\nu = 0.7$}
    \end{subfigure}%
    \begin{subfigure}{0.2\textwidth}
        \centering
        \input{figures/vesicle3d_alpha0.65pro_energy.tex}
        \caption{$\nu = 0.65$}
    \end{subfigure}%
    \begin{subfigure}{0.2\textwidth}
        \centering
        \input{figures/vesicle3d_alpha0.6pro_energy.tex}
        \caption{$\nu = 0.6$}
    \end{subfigure}%
    \caption{
        Evolution of the Willmore energy $\mathcal{E}_{W} = \int_{\Gamma} H^2$ 
        \labelcref{eq:willmore_energy} of the 3D prolate vesicle 
        with the volume- and area-preserving 
        \cref{alg:3dWillmoreFlowWithConservations}. 
    }
    \label{s5:fig:vesicle3dProEnergy}
\end{figure}

\paragraph{}
We can observe for all the -- oblate and prolate -- three-dimensional vesicle
simulations that the computed equilibrium shapes are axisymmetric as expected,
even though our algorithm does not enforce this property. We also note that our
algorithm seems very robust in this case from the energy minimisation point of
view, as it reaches the minimum rather quickly and then preserves it for a long
time without any numerical artifact.

\begin{figure}[!ht]
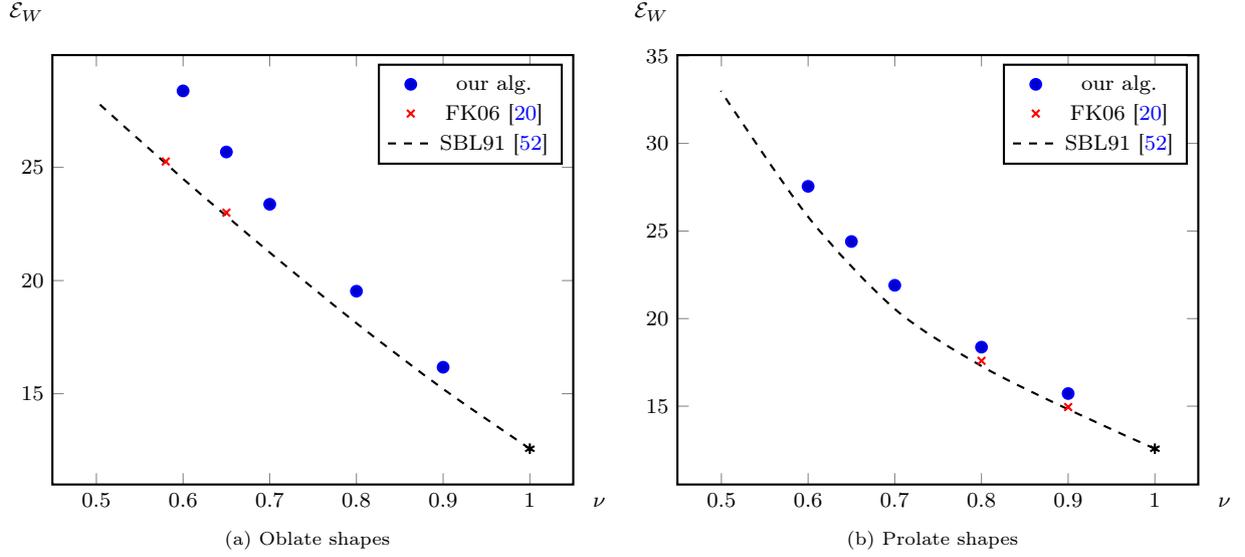

    \centering
    \begin{subfigure}{0.5\textwidth}
        \centering
        \input{figures/vesicle3d_energies_oblate.tex}%
        \caption{Oblate shapes}
    \end{subfigure}%
    \begin{subfigure}{0.5\textwidth}
        \centering
        \input{figures/vesicle3d_energies_prolate.tex}%
        \caption{Prolate shapes}
    \end{subfigure}%
    \caption{
        Willmore energy of the 3D vesicle equilibrium shapes as a function of
        the reduced volume, in the oblate and prolate cases. We compare our
        results to axisymmetric energy minimisation (shooting method) from
        \cite{seifert1991shape} and
        to the ``subdivision thin shell'' vesicle simulations from
        \cite{feng2006finite}.
    }
    \label{s5:fig:vesicle3dEnergiesCmp}
\end{figure}

\begin{figure}[!ht]
    \centering
    \input{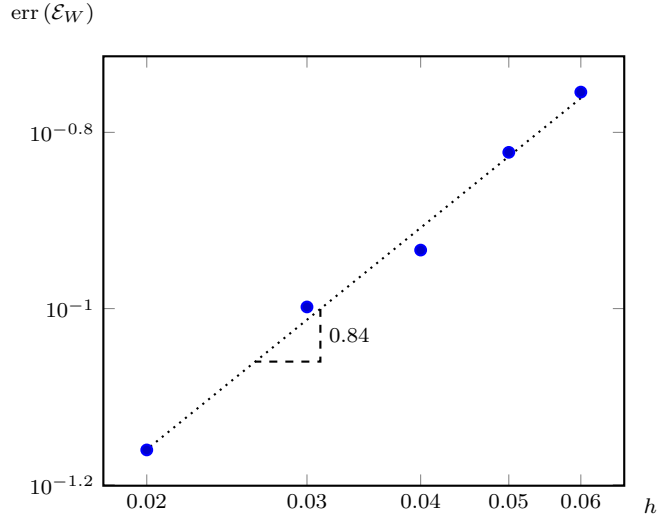}%
    \caption{
        Relative error of the Willmore energy of our three-dimensional vesicle
        equilibrium shape for the $\nu = 0.65$ case as a function of the mesh
        size. The error is computed as $\left( \mathcal{E}_{W} -
        \mathcal{E}_{W}^{th} \right) / \mathcal{E}_{W}^{th}$ with
        $\mathcal{E}_{W}^{th}$ taken from \cite{seifert1991shape}, and the mesh
        sizes range from $h \approx 0.02$ to $h \approx 0.06$. The plot
        is in log-log scale.
    }
    \label{s5:fig:vesicle3dEnergiesConv}
\end{figure}

To compare our results quantitatively, we plot in
\cref{s5:fig:vesicle3dEnergiesCmp} the equilibrium energies of our vesicles
together with the ones obtained in \cite{seifert1991shape} and
\cite{feng2006finite}. We also highlight the numerical convergence of our method
in \cref{s5:fig:vesicle3dEnergiesConv} where we plot the relative error of the
Willmore energy of our final equilibrium shape (for the reduced volume
$\nu=0.65$) as a function of the mesh size. Our results seem in good agreement
with the ones obtained with direct energy minimisation or meshed surface
evolution, despite a small $\sim 10\%$ overestimate for the oblate smallest
reduced volume cases, which seem related to typical finite-element numerical
errors. As shown in \cref{s5:fig:vesicle3dEnergiesConv}, the equilibrium shape
energy converges to the expected value (taken from \cite{seifert1991shape} for
this axisymmetric case) as $h^{0.84}$, which is rather positive for
low-order ($\mathcal{P}^1$) simulations of such high-order effects.

\section{Conclusion}
We propose in this article a numerical method to predict the position of an
interface in 2D or 3D moving according to the gradient of an higher order energy
such as the squared mean curvature (Willmore flow) and subject to constraint of
area and enclosed volume conservation. This problem is of interest to devise
semi-implicit schemes for fluid-structure solvers where the full interaction
between immersed vesicles and fluids are involved. As a first step, the present
work was restricted to the motion of vesicles membranes minimising their mean
curvature at fixed surface area and enclosed volume. The test case considered
corresponds to reach an equilibrium shape starting from an ellipsoid, depending
on its closeness to a sphere or more elongated shape. Even in this delimited
setting, the problem is tricky since the constraints are highly nonlinear and
the energy involved of higher order. Classical approaches through level-set or
phase field methods lead to fourth order PDE to solve, while purely Lagrangian
methods dealing with a surface mesh are not well suited to be included in a
fluid-structure coupling procedure due to the interpolations required between
Eulerian and Lagrangian representations. In this work, we build a method where
only heat equations are solved to compute the right flow of the mean curvature
energy.

Based on diffusion-redistanciation schemes introduced in
\cite{esedog2010diffusion}, our first contribution was to provide a more
intrinsic formulation of these methods, which opens the way to study more easily
other kind of higher order energy. While this study focused on
the mean curvature flows, one could for instance also consider problems where the
Gaussian curvature is involved. We also extended our geometrical flows
algorithms to the evolution of surfaces with conserved area and enclosed volume. 
The originality of our approach relies on a
formula providing explicitly the projection of the unconstrained Willmore flow
on motion conserving area and enclosed volume. This is a big advantage in
comparison to other methods where the area constraint is penalised (and therefore
not exactly fulfilled or leading to stiff problems) or non-linearly enforced (at
a high computational cost).

Our first test case was devoted to the torus, where we observed convergence to
the Clifford torus as expected. Our methodology enjoys numerical convergence and
compares well with existing numerical results to compute equilibrium shapes of
vesicles in dimension $2$. In dimension $3$, while our results are in good
agreement for prolate forms, some discrepancy occurs for oblate vesicle shapes.
This is due to the low order of the finite element method used ($P1$ element),
in contrast with the $C^1$ surface finite element representation of
\cite{feng2006finite}. However we observe numerical convergence in dimension $3$
as well. The use of polynomial discretisation with higher regularity (such as Hermite
elements) would likely
cure that problem without any further modification of our algorithm. One
remarkable feature of our algorithm is to provide consistent results without
any mesh refinement strategy. This is of paramount importance for its possible
use in a full three-dimensional fluid-structure problems, where such remeshing could very
quickly lead, for realistic situations, to intractable computational complexity. 

Another big advantage of our approach is to be very easy to implement: high order
geometric quantities are computed from diffusion equation of a distance
function, which is easily handled by any finite element library.

In a forthcoming work, we will use such a predictor of motion to devise a
semi-implicit scheme for the full fluid-membrane coupling problem.

\section*{Acknowledgements}

Most of the computations presented in this paper were performed using the Froggy
platform of the GRICAD infrastructure (https://gricad.univ-grenoble-alpes.fr),
which is supported by the Rhône-Alpes region (GRANT CPER07\_13 CIRA) and the
Equip@Meso project (reference ANR-10-EQPX-29-01) of the programme
Investissements d’Avenir supervised by the Agence Nationale pour la Recherche,
and the Atlas cluster from the Research Institute in Advanced Mathematics (IRMA
- UMR7501).

\bibliographystyle{plain}
{

\bibliography{refs}
}

\end{document}